\documentclass[12pt]{aptpub}

\usepackage{amsmath,amstext,url}
\usepackage{color}

\numberwithin{equation}{section} 

\makeatletter \@addtoreset{equation}{section}

\makeatletter \@addtoreset{lemma}{section}

\makeatletter \@addtoreset{theorem}{section}

\makeatletter \@addtoreset{proposition}{section}

\makeatletter \@addtoreset{corollary}{section}

\makeatletter \@addtoreset{remark}{section}

\makeatletter \@addtoreset{definition}{section}

\makeatletter \@addtoreset{example}{section}

\begin{document}

\thispagestyle{firstpg}

\vspace*{1.5pc} \noindent \normalsize\textbf{\Large {The jumping properties of Markov branching processes}} \hfill

\vspace{12pt} \hspace*{0.75pc}{\small\textrm{\uppercase{Yanyun Li
}}}\hspace{-2pt}$^{*}$, {\small\textit{ Sun Yat-Sen University}}

\hspace*{0.75pc}{\small\textrm{\uppercase{Junping Li}}}
\hspace{-2pt}$^{**}$, {\small\textit{Central South University }}



\par
\footnote{\hspace*{-0.75pc}$^{*}\,$Postal address:
 School of Mathematics, Sun Yat-Sen University, Guangzhou, 510399, China. E-mail address:
liyy536@mail2.sysu.edu.cn}

\par
\footnote{\hspace*{-0.75pc}$^{**}\,$Postal
address: School of Mathematics and Statistics, Central
South University, Changsha, 410083, China. E-mail address:
jpli@mail.csu.edu.cn}
\par
\par
\renewenvironment{abstract}{%
\vspace{8pt} \vspace{0.1pc} \hspace*{0.25pc}
\begin{minipage}{14cm}
\footnotesize
{\bf Abstract}\\[1ex]
\hspace*{0.5pc}} {\end{minipage}}
\begin{abstract}
  It is well-known that in a Markov branching process, each particle lives for a random long time and gives a random number of new individuals at its death time, independently. It stops when the system has no particle. This paper concentrates on studying the jumping properties of Markov branching process until any time $t\in [0,\infty)$. The joint probability distribution of fixed range jumping numbers of such processes until time $t$ is obtained by using a new method. In particular, the probability distribution of total death number is given for Markov branching processes until time $t\in [0,\infty)$.
\end{abstract}

\vspace*{12pt}
\parbox[b]{26.75pc}{{
}}
\par
\normalsize

\renewcommand{\amsprimary}[1]{
\vspace*{8pt}
\hspace*{2.25pc}
\parbox[b]{20.75pc}{\scriptsize
AMS 2000 Subject Classification: Primary 60J27 Secondary 60J35
     {\uppercase{#1}}}\par\normalsize}
\renewcommand{\ams}[2]{
\vspace*{8pt}
\parbox[b]{24.75pc}{\scriptsize
     MSC 2020 SUBJECT CLASSIFICATION: PRIMARY
          {\uppercase{#1}}\\
 {\uppercase{#2}}}\par\normalsize}

\ams{60J27}

\par
\vspace{5mm}
 \setcounter{section}{1}
 \setcounter{equation}{0}
 \setcounter{theorem}{0}
 \setcounter{lemma}{0}
 \setcounter{corollary}{0}
\noindent {\large \bf 1. Introduction}
\vspace{3mm}
\par
The ordinary Markov branching processes (MBPs) play an important role in the classical field of stochastic processes. The basic property governing the evolution of MBP is the branching property, different individuals act independently when giving birth or death. Each particle in the system lives for a random long time and gives a random number of new individuals at its death time. The system stops when there is no particle in it. It is well-known that $0$ is the absorbing state for Markov branching process. The general discussions of Markov branching processes can be found in Asmussen $\&$ Jagers~\cite{AJ97}, Asmussen $\&$ Hering~\cite{ASHH83}, Athreya $\&$ Ney~\cite{AKN72}, Harris~\cite{Harr63}. Furthermore, some generalized branching systems are studied. For example, Yamazato~\cite{YM75} investigated a branching process with immigration which only occurs at state zero. Chen~\cite{CA2002} and Chen, Li $\&$ Ramesh~\cite{CLR2005} considered general branching processes with or without resurrection. Li $\&$ Chen~\cite{LJCA06} and Li, Chen~$\&$ Pakes~\cite{LCP12} considered branching processes with state-independent immigration. Li $\&$ Liu~\cite{LJL11} considered branching process with migration and immigration. Sevast'yanov~\cite{SBA49} and Vatutin~\cite{VVA74} considered the interacting branching processes. Chen, Li $\&$ Ramesh~\cite{CLR2005} and Chen, Pollett, Zhang $\&$ Li~\cite{CPLZ07} considered weighted Markov branching process. Li, Li $\&$ Chen~\cite{LLC20} investigated the fixed range-jumping properties of weighted Markov branching process.
\par
 However, for a branching system, there are some interesting and important problems remained open. Such as, how many individuals died in the time interval $[0,t)$? What is the $m$-birth number of individuals in the time interval $[0,t)$ (here $m\neq 0$ is fixed)? How many individuals who ever lived in the system (i.e., the total death number) until its extinction? The main purpose of this paper is to consider such problems for Markov branching processes.

\par
To begin our discussion, we first give the definition of Markov branching process. Throughout this paper, let $\mathbb{Z}_+=\{0,1,2,\cdots \}$.
\par
\begin{definition}\label{def1.1}\
A $Q$-matrix $Q=(q_{ij};\ i,j\in \mathbb{Z}_+)$ is called a Markov branching $Q$-matrix (henceforth referred to as MB $Q$-matrix), if
\begin{eqnarray} \label{eq1.1}
 q_{ij}=\begin{cases}ib_{j-i+1},\ & \mbox{if} \ i\geq 1, j\geq i-1, \\
                0 , \ & \mbox{otherwise},
 \end{cases}
 \end{eqnarray}
where
\begin{eqnarray} \label{eq1.2}
           b_j\geq 0 \ (j\neq 1),\  0<-b_1=\sum_{j\neq 1}b_j<\infty.
\end{eqnarray}
\end{definition}
\par
For MB matrix $Q=(q_{ij};\ i,j\in \mathbb{Z}_+)$ defined in Definition~\ref{def1.1}, $b_0$ and $\{b_k; k\geq 2\}$ are the death rate and birth rates of a single individual.
\par
\begin{definition}\label{def1.2}\ A Markov branching process (henceforth referred to as MBP) is a continuous-time Markov chain $\{X(t);\ t\geq 0\}$ with state space $\mathbb{Z}_+$ whose transition function $P(t)=(p_{ij}(t);$ $i,j\in \mathbb{Z}_+)$ satisfies
\begin{eqnarray}\label{eq1.3}
p'_{ij}(t)=\sum_{k=0}^{\infty}p_{ik}(t)q_{kj},\ i\geq0,\ j\geq1,\
t\geq 0,
\end{eqnarray}
where $Q=(q_{ij};\ i,j\in \mathbb{Z}_+)$ is defined in (\ref{eq1.1})-(\ref{eq1.2}).
\end{definition}
\par
By the theory of continuous-time Markov chains, for any given MB $Q$-matrix $Q$, there exists a minimal transition function $P(t)=(p_{ij}(t);$ $i,j\in \mathbb{Z}_+)$, called Feller minimal $Q$-function, satisfying Kolmogorov forward and backward equations. $Q$ is called regular if the Feller minimal $Q$-function $P(t)=(p_{ij}(t);$ $i,j\in \mathbb{Z}_+)$ is honest, i.e., $\sum_{j=0}^{\infty}p_{ij}(t)=1$ for all $i\geq 0$ and $t\geq 0$. Harris~\cite{Harr63} presented the regularity criteria for $Q$ defined in Definition~\ref{def1.1}, i.e., $Q$ is regular if and only if one of the following two conditions holds:
\par
({\bf C}1)\ $B'(1)\leq 0$.
\par
({\bf C}2)\ $B'(1)>0$ and for some (or equivalently, for all) $\varepsilon \in (\rho,1)$,
\begin{eqnarray*}
\int_{\varepsilon}^1\frac{du}{-B(u)}=+\infty,
\end{eqnarray*}
where
\begin{eqnarray}\label{eq1.4}
  B(u)=\sum_{j=0}^{\infty}b_ju^j,\ \ \ u\in [0,1]
\end{eqnarray}
and $\rho$ is the smallest nonnegative root of $B(u)=0$.
\par
The above regularity criteria is named as Harris criteria. We will assume that the process $Q$ is regular in this paper. Therefore, the transition function $P(t)=(p_{ij}(t);$ $i,j\in \mathbb{Z}_+)$ of $\{X(t);t\geq 0\}$ is honest.
\par
\begin{definition}\label{def1.3}
Let $\{X(t);\ t\geq 0\}$ be a Markov branching process defined in Definition~\ref{def1.1} and $i\in \mathbb{Z}_+\setminus\{1\}$. If $X(t_0)-X(t_0-)=i-1$ for some $t_0\in [0,\infty)$, then the jump of $\{X(t); t\geq 0\}$ at time $t_0$ is called an $(i-1)$-range jump. The number of $(i-1)$-range jump of $\{X(t);\ t\geq 0\}$ in time period $[0,t]$ is called $(i-1)$-range jumping number of $\{X(t); t\geq 0\}$ in $[0,t]$.
\end{definition}

\par
\vspace{5mm}
 \setcounter{section}{2}
 \setcounter{equation}{0}
 \setcounter{theorem}{0}
 \setcounter{lemma}{0}
 \setcounter{definition}{0}
 \setcounter{corollary}{0}
\noindent {\large \bf 2. Preliminaries}
 \vspace{3mm}
\par
In this section, we make some preliminaries. Suppose that $D$ is a finite subset of $\mathbb{Z}_+$ with $1\notin D$. Let
\begin{eqnarray*}
  [0,1]^D=\{\emph{\textbf{v}}=(v_k;k\in D); v_k\in [0,1],\ \forall k\in D\}
\end{eqnarray*}
and
\begin{eqnarray*}
  \mathbb{Z}_+^D=\{\emph{\textbf{l}}=(l_k;k\in D); l_k\in \mathbb{Z}_+,\ \forall k\in D\}.
\end{eqnarray*}
For simplicity of notations, in the following, we let $\emph{\textbf{1}}$ denote the vector in $\mathbb{Z}_+^D$ whose components are all $1$ and for $k\in D$, $\emph{\textbf{e}}_k$ denote the vector in $\mathbb{Z}_+^D$ whose $k$'th component is $1$ and others are $0$.
\par
Define
\begin{eqnarray}\label{eq2.1}
B_{D}(u,\emph{\textbf{v}})=\sum_{j\in D}b_j u^j\emph{\textbf{v}}^{\emph{\textbf{e}}_j},\ \ \ \bar{B}_{D}(u)=\sum_{j\in D^c}b_ju^j,\ \ u\in [0,1],\ \emph{\textbf{v}}\in [0,1]^D,
\end{eqnarray}
where $\emph{\textbf{v}}^{\emph{\textbf{l}}}=\prod_{k\in D}v_k^{l_k}$ for $\emph{\textbf{v}}=(v_k;k\in D)$ and $\emph{\textbf{l}}=(l_k;k\in D)$.
\par
It is obvious that $B(u),\ \bar{B}_{D}(u)$ are well defined at least on $[0,1]$, and $B_{D}(u,\emph{\textbf{v}})$ is well defined at least on $[0,1]\times [0,1]^D$.
\par
The following lemma is standard and the proof is omitted.
\par
\begin{lemma} \label{le2.2}
Suppose that $\{f_{\textbf{l}};\ \textbf{l}\in \mathbb{Z}_+^D\}$ is a sequence on $\mathbb{Z}_+^D$, $F(\textbf{v})=\sum_{\textbf{l}\in \mathbb{Z}_+^D}f_{\textbf{l}}\textbf{v}^{\textbf{l}}$ is the generating function of $\{f_{\textbf{l}};\ \textbf{l}\in \mathbb{Z}_+^D\}$. Then for any $j\in \mathbb{Z}_+$,
\begin{eqnarray*}
F^j(\textbf{v})=\sum_{\textbf{l}\in \mathbb{Z}_+^D}f^{*(j)}_{\textbf{l}}\textbf{v}^{\textbf{l}},\ \ \textbf{v}\in [0,1]^D,
\end{eqnarray*}
where  $\textbf{v}^{\textbf{l}}=\prod_{k\in D}v_k^{l_k}$ for $\textbf{v}=(v_k;k\in D)$ and $\textbf{l}=(l_k;k\in D)$, $f^{*(0)}_{\textbf{0}}=1,\ f^{*(0)}_{\textbf{l}}=0\ (\textbf{l}\neq \textbf{0})$ and
\begin{eqnarray*} f^{*(j)}_{\textbf{l}}=\sum_{\textbf{l}^{(1)}+\cdots+\textbf{l}^{(j)}=\textbf{l} }f_{\textbf{l}^{(1)}}\cdots f_{\textbf{l}^{(j)}},\ \ j\geq 1
\end{eqnarray*}
is the $j$'th convolution of $\{f_{\textbf{l}};\ \textbf{l}\in \mathbb{Z}_+^D\}$.
\end{lemma}
\par
The function $\bar{B}_{D}(u)+B_{D}(u,\emph{\textbf{v}})$ will play an important role in our discussion. The following theorem reveals some important properties of this function.
\par
\begin{theorem}\label{th2.1}
{\rm{(i)}}\ For any $\textbf{v}\in[0,1]^{D}$,
\begin{eqnarray}\label{eq2.2}
\bar{B}_{D}(u)+B_{D}(u,\textbf{v})=0
\end{eqnarray}
 possesses at most $2$ roots in $[0,1]$. The minimal nonnegative root of $\bar{B}_{D}(u)+B_{D}(u,\textbf{v})=0$ is denoted by $\rho(\textbf{v})$, then $\rho(\textbf{v})\leq \rho$, where $\rho$ is the minimal nonnegative root of $B(u)=0$.
 \par
 {\rm{(ii)}}\ $\lim_{\textbf{v}\uparrow \textbf{1}}\rho(\textbf{v})=\rho$, where $\textbf{v}\uparrow \textbf{1}$ means $v_k\uparrow 1\ (k\in D)$.
\par
{\rm{(iii)}}\ $\rho(\textbf{v})\in C^{\infty}([0,1)^D)$ and $\rho(\textbf{v})$ can be expanded as a multivariate Taylor series
\begin{eqnarray*}
\rho(\textbf{v})=\sum_{\textbf{l}\in \mathbb{Z}_+^D}\rho_{\textbf{l}}\textbf{v}^{\textbf{l}}.
\end{eqnarray*}
where $\rho_{\textbf{l}}\geq 0, \forall\ \textbf{l}\in \mathbb{Z}_+^D$.
\end{theorem}
\begin{proof}
 Note that $0\leq B_{D}(u,\emph{\textbf{0}})\leq B_{D}(u,\emph{\textbf{v}})\leq B_{D}(u,\emph{\textbf{1}})$, we know that
\begin{eqnarray*}
\bar{B}_{D}(u)+B_{D}(u,\emph{\textbf{v}})\leq B(u).
\end{eqnarray*}
{\rm{(i)}} follows from Li and Chen \cite{Li-Chen08}. Next to prove {\rm{(ii)}}. It is easy to see that $\rho(\emph{\textbf{v}})$ is increasing when $v_k\uparrow $ for all $k\in D$. Therefore, $\lim_{\emph{\textbf{v}}\uparrow \emph{\textbf{1}}}\rho(\emph{\textbf{v}})=\widetilde{\rho}$ exists and $\widetilde{\rho}\leq \rho$. However, if $\widetilde{\rho}<\rho$, substituting $u=\rho(\emph{\textbf{v}})$ in (\ref{eq2.2}) yields
\begin{eqnarray*}
 \bar{B}_{D}(\rho(\emph{\textbf{v}})+B_{D}(\rho(\emph{\textbf{v}}),\emph{\textbf{v}})\equiv 0.
\end{eqnarray*}
Let $\emph{\textbf{v}}\uparrow \emph{\textbf{1}}$ in the above equality, we get
\begin{eqnarray*}
 \bar{B}_{D}(\widetilde{\rho})+B_{D}(\widetilde{\rho},\emph{\textbf{1}})\equiv 0,
\end{eqnarray*}
i.e.,
\begin{eqnarray*}
 B(\widetilde{\rho})=0,
\end{eqnarray*}
which contradicts with minimum property of $\rho$. (ii) is proved.
\par
Now we turn to prove (iii). It follows from Li, Li and Chen~\cite{LLC20} that $\rho(\emph{\textbf{v}})\in C^{\infty}([0,1)^D)$.
\par
Suppose that
\begin{eqnarray*}
\rho(\emph{\textbf{v}})=\sum_{\emph{\textbf{l}}\in \mathbb{Z}_+^D}\rho_{\emph{\textbf{l}}}\emph{\textbf{v}}^{\emph{\textbf{l}}}.
\end{eqnarray*}
\par
Substituting the above expression of $\rho(\emph{\textbf{v}})$ into (\ref{eq2.2}) yields
\begin{eqnarray*}
0& \equiv & \bar{B}_{D}(\rho(\emph{\textbf{v}}))+B_{D}(\rho(\emph{\textbf{v}}),
\emph{\textbf{v}})\\
&=& \sum_{j\in D^c}b_j(\rho(\emph{\textbf{v}}))^j+\sum_{j\in D}b_j(\rho(\emph{\textbf{v}}))^jv_j\\
&=& \sum_{j\in D^c}b_j\sum_{\emph{\textbf{l}}\in \mathbb{Z}_+^D}\rho^{*(j)}_{\emph{\textbf{l}}}
\emph{\textbf{v}}^{\emph{\textbf{l}}}+\sum_{j\in D}b_j\sum_{\emph{\textbf{l}}\in \mathbb{Z}_+^D}\rho^{*(j)}_{\emph{\textbf{l}}}
\emph{\textbf{v}}^{\emph{\textbf{l}}+\emph{\textbf{e}}_j}\\
&=& \sum_{\emph{\textbf{l}}\in \mathbb{Z}_+^D}(\sum_{j\in D^c}b_j\rho^{*(j)}_{\emph{\textbf{l}}})
\emph{\textbf{v}}^{\emph{\textbf{l}}}+\sum_{j\in D}b_j\sum_{\emph{\textbf{l}}\in \mathbb{Z}_+^D}\rho^{*(j)}_{\emph{\textbf{l}}}
\emph{\textbf{v}}^{\emph{\textbf{l}}+\emph{\textbf{e}}_j},
\end{eqnarray*}
where $\emph{\textbf{e}}_j=(\delta_{jk};k\in D)$.
\par
We next prove $\rho_{\emph{\textbf{l}}}\geq 0$ by using mathematical induction respect to $\emph{\textbf{l}}\cdot \emph{\textbf{1}}$. If $\emph{\textbf{l}}\cdot \emph{\textbf{1}}=0$, then $\rho_{\emph{\textbf{0}}}=\rho(\emph{\textbf{0}})\geq 0$ since $\rho(\emph{\textbf{0}})$ is the minimal nonnegative root of $\bar{B}_{ D}(u)+B_{ D}(u,\emph{\textbf{0}})=0$. If $\emph{\textbf{l}}\cdot \emph{\textbf{1}}=1$, i.e., $\emph{\textbf{l}}=\emph{\textbf{e}}_k$ for some $k\in  D$.
Then,
\begin{eqnarray*}
\sum_{j\in  D^c}b_j\rho^{*(j)}_{\emph{\textbf{e}}_k}+b_k\rho^{*(k)}_{\emph{\textbf{0}}}
=0,
\end{eqnarray*}
i.e.,
\begin{eqnarray*}
\sum_{j\in  D^c}b_jj\rho^{j-1}_{\emph{\textbf{0}}}\rho_{\emph{\textbf{e}}_k}
+b_k\rho^{k}_{\emph{\textbf{0}}}
=0.
\end{eqnarray*}
Hence
\begin{eqnarray*}
\rho_{\emph{\textbf{e}}_k}
=-\frac{b_k\rho^{k}_{\emph{\textbf{0}}}}{\bar{B}'_{ D}
(\rho_{\emph{\textbf{0}}})}\geq 0,\ \ k\in  D,
\end{eqnarray*}
since $\bar{B}'_{ D}(\rho_{\emph{\textbf{0}}})<0$.
\par
Assume $\rho_{\emph{\textbf{l}}}\geq 0$ for $\emph{\textbf{l}}$ satisfying $\emph{\textbf{l}}\cdot \emph{\textbf{1}}\leq m$, then for any $\bar{\emph{\textbf{l}}}\in \mathbb{Z}_+^D$ satisfying $\bar{\emph{\textbf{l}}}\cdot \emph{\textbf{1}}=m+1$, there exists $\emph{\textbf{l}}\in \mathbb{Z}_+^D$ and $k\in  D$ such that $\bar{\emph{\textbf{l}}}=\emph{\textbf{l}}+\emph{\textbf{e}}_k$ and $\emph{\textbf{l}}\cdot \emph{\textbf{1}}\leq m$, therefore,
\begin{eqnarray*}
\sum_{j\in  D^c}b_j\rho^{*(j)}_{\emph{\textbf{l}}+\emph{\textbf{e}}_k}
+b_k\rho^{*(k)}_{\emph{\textbf{l}}}=0,
\end{eqnarray*}
i.e.,
\begin{eqnarray*}
\sum_{j\in  D^c}b_jj\rho^{j-1}_{\emph{\textbf{0}}}\rho_{\emph{\textbf{l}}
+\emph{\textbf{e}}_k}
+\sum_{j\in  D^c\setminus\{1\}}b_j\sum_{\emph{\textbf{l}}^{(1)}+\cdots +\emph{\textbf{l}}^{(j)}=\emph{\textbf{l}}+\emph{\textbf{e}}_k,\ \emph{\textbf{l}}^{(1)}\cdot\emph{\textbf{1}},\cdots ,\emph{\textbf{l}}^{(j)}\cdot \emph{\textbf{1}}\leq m}\rho_{\emph{\textbf{l}}^{(1)}}
\cdots \rho_{\emph{\textbf{l}}^{(j)}}+b_k\rho^{*(k)}_{\emph{\textbf{l}}}=0.
\end{eqnarray*}
Hence
\begin{eqnarray*}
\rho_{\bar{\emph{\textbf{l}}}}=\rho_{\emph{\textbf{l}}+\emph{\textbf{e}}_k}
=-\frac{\sum_{j\in  D^c\setminus\{1\}}b_j\sum_{\emph{\textbf{l}}^{(1)}+\cdots +\emph{\textbf{l}}^{(j)}=\emph{\textbf{l}}+\emph{\textbf{e}}_k,\ \emph{\textbf{l}}^{(1)}\cdot\emph{\textbf{1}},\cdots ,\emph{\textbf{l}}^{(j)}\cdot \emph{\textbf{1}}\leq m}\rho_{\emph{\textbf{l}}^{(1)}}
\cdots \rho_{\emph{\textbf{l}}^{(j)}}+b_k\rho^{*(k)}_{\emph{\textbf{l}}}}{\bar{B}'_{
 D}(\rho_{\emph{\textbf{0}}})}\geq 0,
\end{eqnarray*}
since $\bar{B}'_{ D}(\rho_{\emph{\textbf{0}}})<0$.
By mathematical induction, we know that $\rho_{\emph{\textbf{l}}}\geq 0, \forall\ \emph{\textbf{l}}\in \mathbb{Z}_+^D$. The proof is complete. \hfill $\Box$
\end{proof}

\par
\vspace{5mm}
 \setcounter{section}{3}
 \setcounter{equation}{0}
 \setcounter{theorem}{0}
 \setcounter{lemma}{0}
 \setcounter{definition}{0}
 \setcounter{corollary}{0}
\noindent {\large \bf 3. Down/up jumping property}
 \vspace{3mm}
\par
We now turn to consider the down/up jumping properties of $\{X(t);t\geq 0\}$.
\par
As in the previous section, let $ D\subset \mathbb{Z}_+$ be a finite subset with $1\notin  D$ and $d=| D|$ denote the number of elements in $D$. We also assume that $b_k>0$ for all $k\in D$ since $\{X(t);t\geq 0\}$ has no $(k-1)$-range jump if $b_k=0$. For simplicity of notation, denote
$$
D-1=\{k-1; k\in D\}.
$$
\par
 The main purpose of this paper is to consider the $( D-1)$-range jumping numbers of $\{X(t);t\geq 0\}$. However, the process $\{X(t);t\geq 0\}$ itself can not reveal such jumping numbers directly. Therefore, we need to find a new method to discuss the property of $(D-1)$-range jumping numbers. For this purpose, we construct a new $Q$-matrix $\widetilde{Q}=(\tilde{q}_{(i,\emph{\textbf{m}}), (j,\emph{\textbf{l}})};\ (i,\emph{\textbf{m}}),$ $ (j,\emph{\textbf{l}})\in \mathbb{Z}_+\times\mathbb{Z}_+^D)$, where
\begin{eqnarray}\label{eq3.1}
 \tilde{q}_{(i,\emph{\textbf{m}}), (j,\emph{\textbf{l}})}=\begin{cases}ib_{j-i+1},\ & if \ i\geq 1, j-i+1\in  D^c, \emph{\textbf{l}}=\emph{\textbf{m}} \\
                   ib_{j-i+1}, \ & if \ i\geq 1, j-i+1\in  D, \emph{\textbf{l}}=\emph{\textbf{m}}+\emph{\textbf{e}}_{j-i+1}\\
                0 , \ & otherwise
                     \end{cases}
\end{eqnarray}
and $\{b_k;k\geq 0\}$ is given in (\ref{eq1.3}).

\par
Let $\widetilde{P}(t)=(\tilde{p}_{(i,\emph{\textbf{m}}), (j,\emph{\textbf{l}})}(t); (i,\emph{\textbf{m}}), (j,\emph{\textbf{l}})\in \mathbb{Z}_+\times\mathbb{Z}_+^D)$ denote the Feller minimal $\widetilde{Q}$-function and
\begin{eqnarray*}
F_{i,\emph{\textbf{m}}}(t,u,\emph{\textbf{v}})=
\sum_{(j,\emph{\textbf{l}})\in \mathbb{Z}_+\times\mathbb{Z}_+^D}\tilde{p}_{(i,\emph{\textbf{m}}),
 (j,\emph{\textbf{l}})}(t)u^j\emph{\textbf{v}}^{\emph{\textbf{l}}},\ \ (u,\emph{\textbf{v}})\in [0,1]\times [0,1]^D,
\end{eqnarray*}
where $\emph{\textbf{v}}^{\emph{\textbf{l}}}=\prod_{k\in  D}v_k^{l_k}$ for $\emph{\textbf{v}}=(v_k;k\in D)$ and $\emph{\textbf{l}}=(l_k;k\in D)$.
\par
\begin{lemma}\label{le3.1a} Let $\widetilde{Q}$ be defined in $(\ref{eq3.1})$ and $\widetilde{P}(t)=(\tilde{p}_{(i,\textbf{m}), (j,\textbf{l})}(t); (i,\textbf{m}), (j,\textbf{l})\in \mathbb{Z}_+\times\mathbb{Z}_+^D)$ be the Feller minimal $\widetilde{Q}$-function. Then
\par
{\rm{(i)}}\ for any $(i,\textbf{m})\in \mathbb{Z}_+\times\mathbb{Z}_+^D$,
\begin{eqnarray}\label{eq3.2}
\frac{\partial F_{i,\textbf{m}}(t,u,\textbf{v})}{\partial t}=[\bar{B}_{D}(u)+B_{D}
 (u,\textbf{v})]\cdot \frac{\partial F_{i,\textbf{m}}(t,u,\textbf{v})}{\partial u},\ \ (u,\textbf{v})\in [0,1]\times[0,1]^D,
\end{eqnarray}
where $\bar{B}_{D}(u),\ B_{D}(u,\textbf{v})$ are defined in $(\ref{eq2.1})$. Moreover,
\begin{eqnarray}\label{eq3.3}
 F_{i,\textbf{m}}(t,u,\textbf{v})-u^i= [\bar{B}_{D}(u)+B_{D}(u,\textbf{v})]\cdot\int_0^t\frac{\partial F_{i,\textbf{m}}(s,u,\textbf{v})}{\partial u}ds.
\end{eqnarray}
\par
{\rm{(ii)}}\ $\widetilde{Q}$ is regular if and only if $Q$ is regular, i.e., $({\bf C}1)$ or $({\bf C}2)$ holds.
\end{lemma}
\par
\begin{proof}
It follows from Kolmogorov forward equations that
 \begin{eqnarray*}
 &&\sum_{(j,\emph{\textbf{l}})\in \mathbb{Z}_+\times\mathbb{Z}_+^D}
 \tilde{p}'_{(i,\emph{\textbf{m}}),(j,\emph{\textbf{l}})}(t)
 u^j\emph{\textbf{v}}^{\emph{\textbf{l}}}\\
 &=&\sum_{(k,\emph{\textbf{r}})\in \mathbb{Z}_+\times \mathbb{Z}_+^D}\tilde{p}_{(i,\emph{\textbf{m}}),
 (k,\emph{\textbf{r}})}(t)\cdot \sum_{(j,\emph{\textbf{l}})\in \mathbb{Z}_+\times\mathbb{Z}_+^D}
 \tilde{q}_{(k,\emph{\textbf{r}}),(j,\emph{\textbf{l}})}
 u^j\emph{\textbf{v}}^{\emph{\textbf{l}}}\\
 &=&\sum_{(k,\emph{\textbf{r}})\in \mathbb{Z}_+\times \mathbb{Z}_+^D}\tilde{p}_{(i,\emph{\textbf{m}}),
 (k,\emph{\textbf{r}})}(t)\cdot [\sum_{j:j-k+1\in D}
 kb_{j-k+1}u^j\emph{\textbf{v}}^{\emph{\textbf{r}}
 +\emph{\textbf{e}}_{j-k+1}}+\sum_{j:j-k+1\notin D}kb_{j-k+1}u^j\emph{\textbf{v}}^{\emph{\textbf{r}}}]\\
 &=&\sum_{(k,\emph{\textbf{r}})\in \mathbb{Z}_+\times \mathbb{Z}_+^D}\tilde{p}_{(i,\emph{\textbf{m}}),
 (k,\emph{\textbf{r}})}(t)\cdot ku^{k-1}\emph{\textbf{v}}^{\emph{\textbf{r}}}
 \cdot[\sum_{j\in D}
 b_{j}u^j\emph{\textbf{v}}^{\emph{\textbf{e}}_j}
 +\sum_{j\in D^c}b_ju^j].
 \end{eqnarray*}
(i) is proved. We now prove (ii).
\par
If $B'(1)\leq 0$, then $\rho=1$. By Harris criteria, $Q$ is regular. On the other hand, let $u=\rho(\emph{\textbf{v}})$ in (\ref{eq3.3}), we know that
\begin{eqnarray*}
 F_{i,\emph{\textbf{m}}}(t,\rho(\emph{\textbf{v}}),\emph{\textbf{v}})
 -\rho^i(\emph{\textbf{v}})=0.
\end{eqnarray*}
Let $\emph{\textbf{v}}\uparrow \emph{\textbf{1}}$ in the above equality and use monotone convergence theorem, we get $F_{i,\emph{\textbf{m}}}(t,1,\emph{\textbf{1}})=1$, i.e., $\widetilde{Q}$ is regular. Therefore, we only need to consider the case $B'(1)>0$.
\par
Let
\begin{eqnarray*}
 \tilde{\phi}_{(i,\emph{\textbf{m}}),(j,\emph{\textbf{l}})}(\lambda)=\int_0^{\infty}
 e^{\lambda t}\cdot \tilde{p}_{(i,\emph{\textbf{m}}),(j,\emph{\textbf{l}})}(t)dt
\end{eqnarray*}
be the Laplace transform of $\tilde{p}_{(i,\emph{\textbf{m}}),(j,\emph{\textbf{l}})}(t)$. It follows from (\ref{eq3.2}) that
\begin{eqnarray}\label{eq3.4}
\lambda \tilde{\Phi}_{i,\emph{\textbf{m}}}(\lambda,u,\emph{\textbf{v}})-u^i
\emph{\textbf{v}}^{\emph{\textbf{m}}}=[\bar{B}_{D}(u)+B_{D}
 (u,\emph{\textbf{v}})]\cdot \frac{\partial \tilde{\Phi}_{i,\emph{\textbf{m}}}(\lambda,u,\emph{\textbf{v}})}{\partial u},
\end{eqnarray}
where $\tilde{\Phi}_{i,\emph{\textbf{m}}}(\lambda,u,\emph{\textbf{v}})=\sum_{(j,\emph{\textbf{l}})\in \mathbb{Z}_+\times\mathbb{Z}_+^D}\tilde{\phi}_{(i,\emph{\textbf{m}}),
 (j,\emph{\textbf{l}})}(\lambda)u^j\emph{\textbf{v}}
 ^{\emph{\textbf{l}}}$.
 \par
 Suppose that $Q$ is regular but $\widetilde{Q}$ is not regular.
Then there exist some $(i,\emph{\textbf{m}})\in \mathbb{Z}_+\times\mathbb{Z}_+^D$ and some $\lambda>0$ such that $1-\lambda \tilde{\Phi}_{i,\emph{\textbf{m}}}(\lambda,1,\emph{\textbf{1}})
=\delta(i,\emph{\textbf{m}},\lambda)>0$ and hence, there exists $\varepsilon in (\rho,1)$ such that $u^i-\lambda \tilde{\Phi}_{i,\emph{\textbf{m}}}(\lambda,u,\emph{\textbf{1}})
\geq \frac{\delta(i,\emph{\textbf{m}},\lambda)}{2}$ for all $u\in (\varepsilon, 1]$. It follows from (\ref{eq3.4}) that
\begin{eqnarray*}
&&\tilde{\Phi}_{i,\emph{\textbf{m}}}(\lambda,1,\emph{\textbf{1}})
-\tilde{\Phi}_{i,\emph{\textbf{m}}}(\lambda,\varepsilon,\emph{\textbf{1}})\\
&=&\int_{\varepsilon}^1\frac{\lambda \tilde{\Phi}_{i,\emph{\textbf{m}}}(\lambda,u,\emph{\textbf{1}})-u^i}
{B(u)}du\\
 &\geq &\frac{\delta(i,\emph{\textbf{m}},\lambda)}{2}\cdot
 \int_{\varepsilon}^1\frac{du}
{-B(u)}\\
&=&+\infty
\end{eqnarray*}
which is contradicts with $\lambda \tilde{\Phi}_{i,\emph{\textbf{m}}}(\lambda,1,\emph{\textbf{1}})\leq 1$. Therefore, $\widetilde{Q}$ is also regular.
\par
Conversely, suppose that $\widetilde{Q}$ is regular. By the theory of Markov chains (see Anderson~\cite{And91}), the $\widetilde{Q}$-function $(\tilde{p}_{(i,\emph{\textbf{m}}),
 (j,\emph{\textbf{l}})}(t);(i,\emph{\textbf{m}}),
 (j,\emph{\textbf{l}})\in \mathbb{Z}_+\times\mathbb{Z}_+^D)$ can be obtained as follows.
 \begin{eqnarray*}
 &&\ \ \ \ \tilde{p}^{(n)}_{(i,\emph{\textbf{m}}),
 (j,\emph{\textbf{l}})}(t)\\
 &&=
 \begin{cases}
 \delta_{(i,\emph{\textbf{m}}),
 (j,\emph{\textbf{l}})}e^{-\tilde{q}_{(i,\emph{\textbf{m}})}t},\ \ & if\ n=0\\
 \tilde{p}^{(0)}_{(i,\emph{\textbf{m}}),
 (j,\emph{\textbf{l}})}(t)+\int_0^te^{-\tilde{q}_{(i,\emph{\textbf{m}})}s}
 \sum_{(k,\emph{\textbf{r}})\neq(i,\emph{\textbf{m}})}
 \tilde{q}_{(i,\emph{\textbf{m}}),
 (k,\emph{\textbf{r}})}\cdot\tilde{p}^{(n-1)}_{(k,\emph{\textbf{r}}),
 (j,\emph{\textbf{l}})}(t-s)ds,\ \ & if\ n\geq 1\\
 \end{cases}\\
 && \tilde{p}_{(i,\emph{\textbf{m}}),
 (j,\emph{\textbf{l}})}(t)=\lim_{n\rightarrow \infty}\tilde{p}^{(n)}_{(i,\emph{\textbf{m}}),
 (j,\emph{\textbf{l}})}(t),\ \ (i,\emph{\textbf{m}}),
 (j,\emph{\textbf{l}})\in \mathbb{Z}_+\times\mathbb{Z}_+^D,
 \end{eqnarray*}
 where $\tilde{q}_{(i,\emph{\textbf{m}})}
 =-\tilde{q}_{(i,\emph{\textbf{m}}),(i,\emph{\textbf{m}})}$. For any $i,j\in \mathbb{Z}_+,\ \emph{\textbf{m}}\in \mathbb{Z}_+^D$, let
 $$
 f^{(n)}_{i,j}(t,\emph{\textbf{m}})=\sum_{\emph{\textbf{l}}\in \mathbb{Z}_+^D}\tilde{p}^{(n)}_{(i,\emph{\textbf{m}}),
 (j,\emph{\textbf{l}})}(t),\ \ n\geq 0.
 $$
 By the definition of $\widetilde{Q}$, we can see that $\tilde{q}_{(i,\emph{\textbf{m}})}=-ib_1$ and hence $f^{(0)}_{ij}(t,\emph{\textbf{m}})=\delta_{ij}e^{ib_1t}$ which is independent of $\emph{\textbf{m}}$. Thus, we can rewrite $f^{(0)}_{ij}(t,\emph{\textbf{m}})$ as $f^{(0)}_{ij}(t)$. Suppose that $f^{(n-1)}_{ij}(t,\emph{\textbf{m}})$ is dependent of $\emph{\textbf{m}}$, rewritten as $f^{(n-1)}_{ij}(t)$. Then,
 \begin{eqnarray*}
 &&f^{(n)}_{ij}(t,\emph{\textbf{m}})\\
 &=& f^{(0)}_{ij}(t)+\int_0^te^{ib_1s}
 \sum_{(k,\emph{\textbf{r}})\neq(i,\emph{\textbf{m}})}
 \tilde{q}_{(i,\emph{\textbf{m}}),
 (k,\emph{\textbf{r}})}\cdot f^{(n-1)}_{kj}(t-s)ds\\
 &=& f^{(0)}_{ij}(t)+\int_0^te^{ib_1s}
 [\sum_{\emph{\textbf{r}}\neq\emph{\textbf{m}}}
 \tilde{q}_{(i,\emph{\textbf{m}}),
 (i,\emph{\textbf{r}})}\cdot f^{(n-1)}_{ij}(t-s)+\sum_{k\neq i}\sum_{\emph{\textbf{r}}\in \mathbb{Z}_+^D}\tilde{q}_{(i,\emph{\textbf{m}}),
 (k,\emph{\textbf{r}})}\cdot f^{(n-1)}_{kj}(t-s)]ds\\
 &=& f^{(0)}_{ij}(t)+\int_0^te^{ib_1s}
 \sum_{k\geq i-1}ib_{k-i+1}\cdot f^{(n-1)}_{kj}(t-s)ds\\
 \end{eqnarray*}
 which implies that $f^{(n)}_{ij}(t,\emph{\textbf{m}})$ is dependent of $\emph{\textbf{m}}$. By the mathematical induction, we know that for any $n\geq 0$, $f^{(n)}_{ij}(t,\emph{\textbf{m}})$ is dependent of $\emph{\textbf{m}}$ and can be rewritten as $f^{(n)}_{ij}(t)$. In fact, we have proved that $f^{(n)}_{ij}(t)$ satisfies
  \begin{eqnarray*}
 &&\ \ \ \ f^{(n)}_{ij}(t)\\
 &&=
 \begin{cases}
 \delta_{ij}e^{ib_1t},\ \ & if\ n=0\\
 f^{(0)}_{ij}(t)+\int_0^te^{ib_1s}
 \sum_{k\geq i-1}
 ib_{k-i+1}\cdot f^{(n-1)}_{kj}(t-s)ds,\ \ & if\ n\geq 1.\\
 \end{cases}
 \end{eqnarray*}
By the theory of Markov chains, we know that $f_{ij}(t)=\lim_{n\rightarrow \infty}f^{(n)}_{ij}(t)$ exists and is the Feller minimal $Q$-function. Furthermore,
\begin{eqnarray*}
\sum_{j=0}^{\infty}f_{ij}(t)
   &=&
   \lim_{n\rightarrow\infty}\sum_{j=0}^{\infty}f^{(n)}_{ij}(t)\\
   &=&
   \lim_{n\rightarrow\infty}
   \sum_{j=0}^{\infty}f^{(n)}_{ij}(t,\emph{\textbf{m}})\\
   &=&
   \lim_{n\rightarrow\infty}
   \sum_{j=0}^{\infty}\sum_{\emph{\textbf{l}}\in \mathbb{Z}_+^D}\tilde{p}^{(n)}_{(i,\emph{\textbf{m}}),
   (j,\emph{\textbf{l}})}(t)\\
   &=&
   \sum_{j=0}^{\infty}\sum_{\emph{\textbf{l}}\in \mathbb{Z}_+^D}\tilde{p}_{(i,\emph{\textbf{m}}),
   (j,\emph{\textbf{l}})}(t)=1.\\
\end{eqnarray*}
 Therefore, $Q$ is regular. The proof is complete.\hfill $\Box$
\end{proof}
\par
Since we have assumed that $Q$ is regular, by the definition of $\widetilde{Q}$ and Lemma~\ref{le3.1a}, we see that $\widetilde{Q}$ determines a unique $\widetilde{Q}$-process $\{(\tilde{X}(t),\emph{\textbf{Y}}(t)); t\geq 0\}$, where $\emph{\textbf{Y}}(t)=(Y_k(t); k\in D)$ counts the $(D-1)$-jumping number of $\{\tilde{X}(t);t\geq 0\}$. It follows from the proof of Lemma~\ref{le3.1a} that $\{\tilde{X}(t);t\geq 0\}$ is the MBP with generator $Q$ and hence has the same distribution as $\{X(t);t\geq 0\}$. Therefore, we still use $\{X(t);t\geq 0\}$ to denote $\{\tilde{X}(t);t\geq 0\}$ in the following. i.e., $\{(X(t),\emph{\textbf{Y}}(t)); t\geq 0\}$ is the $\widetilde{Q}$-process, where $\{X(t);t\geq 0\}$ is the MBP defined in Definition~\ref{def1.2} and $\emph{\textbf{Y}}(t)=(Y_k(t); k\in D)$ counts the $(D-1)$-jumping number of $\{X(t);t\geq 0\}$. In particular,
\par
(i)\ if $ D=\{0\}$ then $Y_0(t)$ counts the down jumping number (i.e., the death number) of $\{X(t): t\geq 0\}$ until time $t$.
\par
(ii)\ If $ D=\{i\}\ (i\geq 2)$, then $Y_i(t)$ counts the $(i-1)$-range up jumping number of $\{X(t): t\geq 0\}$ until time $t$.
\par
(iii)\ If $ D=\{0,i\}\ (i\geq 2)$, then $\emph{\textbf{Y}}(t)=(Y_0(t),Y_i(t))$ counts the death number and the $(i-1)$-range up jumping number of $\{X(t): t\geq 0\}$ until time $t$.

\par
\begin{lemma}\label{le3.1}
For $\widetilde{P}(t)$, we have that for any $(i,\textbf{m})\in \mathbb{Z}_+\times\mathbb{Z}_+^D$,
\begin{eqnarray}\label{eq3.5}
F_{i,\textbf{m}}(t,u,\textbf{v})=[F_{1,\textbf{0}}(t,u,\textbf{v})]^i\cdot \textbf{v}^{\textbf{m}},\ \ \ (u,\textbf{v})\in [0,1]\times[0,1]^D,
\end{eqnarray}
where $\textbf{v}^{\textbf{m}}=\prod_{k\in D}v_k^{m_k}$ for $\textbf{v}=(v_k;k\in D)$ and $\textbf{m}=(m_k;k\in D)$.
\end{lemma}
\begin{proof}
\par
For any $i\geq 0,\ \emph{\textbf{m}}\in \mathbb{Z}_+^D$, let $X_k(t)$ be the offspring number at time $t$ of $k$th individual, $\emph{\textbf{Y}}^{(k)}(t)$ be the $( D-1)$ jumping number of $\{X_k(t);t\geq 0\}$. Then, $\{(X_k(t),\emph{\textbf{Y}}^{(k)}(t));\ k=1,\cdots, i\}$ are independent and identically distributed with the same distribution as $(X(t),\emph{\textbf{Y}}(t))$ starting at $(X(0),\emph{\textbf{Y}}(0))=(1,\emph{\textbf{0}})$. Therefore,
\begin{eqnarray*}
&& E[u^{X(t)}\cdot \emph{\textbf{v}}^{\emph{\textbf{Y}}(t)}|(X(0),\emph{\textbf{Y}}(0))
=(i,\emph{\textbf{m}})]\\
&=& E[u^{\sum_{k=1}^iX_k(t)}\cdot \emph{\textbf{v}}^{\emph{\textbf{m}}+\sum_{k=1}^i\emph{\textbf{Y}}^{(k)}(t)}]\\
 &=& E[\prod_{k=1}^iu^{X_k(t)}\cdot \emph{\textbf{v}}^{\emph{\textbf{Y}}^{(k)}(t)}]
 \cdot\emph{\textbf{v}}^{\emph{\textbf{m}}}\\
 &=& \left(\prod_{k=1}^iE[u^{X_k(t)}\cdot \emph{\textbf{v}}^{\emph{\textbf{Y}}^{(k)}(t)}]\right)
 \cdot\emph{\textbf{v}}^{\emph{\textbf{m}}}\\
 &=& \left(E[u^{X(t)}\cdot \emph{\textbf{v}}^{\emph{\textbf{Y}}(t)}|(X(0),\emph{\textbf{Y}}(0))
=(1,\emph{\textbf{0}})]\right)^i
 \cdot\emph{\textbf{v}}^{\emph{\textbf{m}}}.\\
\end{eqnarray*}
which achieves (\ref{eq3.5}). The proof is complete.
\hfill $\Box$
\end{proof}
\par
\begin{lemma}\label{le3.2} Suppose that $\textbf{v}\in [0,1)^D$ and $u\in [0,1]$.
\par
{\rm{(i)}}\ The differential equation
\begin{eqnarray}\label{eq3.6}
\begin{cases}\frac{\partial y}{\partial t}=B_{ D}(y,\textbf{v})+\bar{B}_{ D}(y)\\
y|_{t=0}=u
\end{cases}
\end{eqnarray}
has unique solution $y(t;u,\textbf{v})=G(t,u,\textbf{v})$.
\par
{\rm{(ii)}}\ If $u\in [0,\rho(\textbf{v}))$, then $y(t;u,\textbf{v})=G(t,u,\textbf{v})$ is increasing to $\rho(\textbf{v})$ as $t\uparrow \infty$. If $u\in (\rho(\textbf{v}),1]$, then $y(t;u,\textbf{v})=G(t,u,\textbf{v})$ is decreasing to $\rho(\textbf{v})$ as $t\uparrow \infty$. If $u=\rho(\textbf{v})$, then $G(t,u,\textbf{v})\equiv\rho(\textbf{v})$.
\end{lemma}
\begin{proof}
We first prove the existence of solution to (\ref{eq3.6}). Denote $H(x)=B_{ D}(x,\emph{\textbf{v}})+\bar{B}_{ D}(x)-b_1x$. Take $y_0(t;u,\\emph{\textbf{v}})\equiv u$ and let
$$
y_n(t;u,\emph{\textbf{v}})=e^{b_1t}\cdot[u+\int_0^te^{-b_1s}H(y_{n-1}(s;u,\emph{\textbf{v}}))ds],\ \ n\geq 1.
$$
\par
(a)\ If $u\in (\rho(\emph{\textbf{v}}),1]$, it can be proved that $y_n(t;u,\emph{\textbf{v}})> \rho(\emph{\textbf{v}})$ and $y_n(t;u,\emph{\textbf{v}})\leq y_{n-1}(t;u,\emph{\textbf{v}})\ (n\geq 1)$. Indeed, obviously, $y_0(t)\equiv u>\rho(\emph{\textbf{v}})$. Assume $y_n(t;u,\emph{\textbf{v}})>\rho(\emph{\textbf{v}})$, then
\begin{eqnarray*}
y_{n+1}(t;u,\emph{\textbf{v}})&=&e^{b_1t}\cdot[u+\int_0^te^{-b_1s}
H(y_n(s;u,\emph{\textbf{v}}))ds]\\
&> & e^{b_1t}\cdot[u+\int_0^te^{-b_1s}H(\rho(\emph{\textbf{v}}))ds]\\
&= & e^{b_1t}\cdot[u-b_1\rho(\emph{\textbf{v}})\int_0^te^{-b_1s}ds]\\
&= & e^{b_1t}\cdot[u+\rho(\emph{\textbf{v}})e^{-b_1t}-\rho(\emph{\textbf{v}})]\\
&> & \rho(\emph{\textbf{v}}).
\end{eqnarray*}
On the other hand, $y_1(t;u,\emph{\textbf{v}})=e^{b_1t}\cdot[u+\int_0^te^{-b_1s}H(u)ds]< e^{b_1t}\cdot[u-b_1u\int_0^te^{-b_1s}ds]=u=y_0(t;u,\emph{\textbf{v}})$. Assume $y_n(t;u,\emph{\textbf{v}})\leq y_{n-1}(t;u,\emph{\textbf{v}})$, then,
\begin{eqnarray*}
y_{n+1}(t;u,\emph{\textbf{v}})&=&e^{b_1t}\cdot[u+\int_0^te^{-b_1s}
H(y_n(s;u,\emph{\textbf{v}}
))ds]\\
&\leq & e^{b_1t}\cdot[u+\int_0^te^{-b_1s}H(y_{n-1}(s;u,\emph{\textbf{v}}))ds]\\
&= & y_n(t;u,\emph{\textbf{v}}).
\end{eqnarray*}
Therefore, it follows from monotone convergence theorem that
$G(t,u,\emph{\textbf{v}})=\lim_{n\rightarrow \infty}y_n(t;u,\emph{\textbf{v}})$ exists and satisfies
$$
G(t,u,\emph{\textbf{v}})=e^{b_1t}\cdot[u+\int_0^te^{-b_1s}H(G(s,u,\emph{\textbf{v}}))
ds].
$$
Hence, $y(t;u,\emph{\textbf{v}})=G(t,u,\emph{\textbf{v}})$ is a solution of (\ref{eq3.6}). Since $B_{D}(y,\emph{\textbf{v}})+\bar{B}_{ D}(y)<0$ for all $y\in (\rho(\emph{\textbf{v}}),1]$, we know that $G(t,u,\emph{\textbf{v}})$ is decreasing and it is easy to see that $\lim_{t\rightarrow \infty}G(t,u,\emph{\textbf{v}})=\rho(\emph{\textbf{v}})$.
\par
(b)\ If $u\in [0,\rho(\emph{\textbf{v}}))$, it can be proved that $y_n(t;u,\emph{\textbf{v}})< \rho(\emph{\textbf{v}})$ and $y_n(t;u,\emph{\textbf{v}})\geq y_{n-1}(t;u,\emph{\textbf{v}})\ (n\geq 1)$. Indeed, obviously, $y_0(t;u,\emph{\textbf{v}})\equiv u<\rho(\emph{\textbf{v}})$. Assume $y_n(t;u,\emph{\textbf{v}})<\rho(\emph{\textbf{v}})$, similar as in (a),
\begin{eqnarray*}
y_{n+1}(t;u,\emph{\textbf{v}})&<& e^{b_1t}\cdot[u+\int_0^te^{-b_1s}H(\rho(\emph{\textbf{v}}))ds]\\
&<& \rho(\emph{\textbf{v}}).
\end{eqnarray*}
On the other hand, $y_1(t;u,\emph{\textbf{v}})=e^{b_1t}\cdot[u+\int_0^te^{-b_1s}H(u)ds]> e^{b_1t}\cdot[u-b_1u\int_0^te^{-b_1s}ds]=u=y_0(t;u,\emph{\textbf{v}})$. Assume $y_n(t;u,\emph{\textbf{v}})\geq y_{n-1}(t;u,\emph{\textbf{v}})$, then,
\begin{eqnarray*}
y_{n+1}(t;u,\emph{\textbf{v}})&=&e^{b_1t}\cdot[u+\int_0^te^{-b_1s}
H(y_n(s;u,\emph{\textbf{v}}))ds]\\
&\geq & e^{b_1t}\cdot[u+\int_0^te^{-b_1s}H(y_{n-1}(s;u,\emph{\textbf{v}}))ds]\\
&= & y_n(t).
\end{eqnarray*}
Therefore, it follows from monotone convergence theorem that
$G(t,u,\emph{\textbf{v}})=\lim_{n\rightarrow \infty}y_n(t;u,\emph{\textbf{v}})$ exists and satisfies
$$
G(t,u,\emph{\textbf{v}})=e^{b_1t}\cdot[u+\int_0^te^{-b_1s}H(G(s,u,\emph{\textbf{v}}))
ds].
$$
Hence, $y(t;u,\emph{\textbf{v}})=G(t,u,\emph{\textbf{v}})$ is a solution of (\ref{eq3.6}). Since $B_{ D}(y,\emph{\textbf{v}})+\bar{B}_{ D}(y)>0$ for all $y\in [0,\rho(\emph{\textbf{v}}))$, we know that $G(t,u,\emph{\textbf{v}})$ is increasing and it is easy to see that $\lim_{t\rightarrow\infty}G(t,u,\emph{\textbf{v}})=\rho(\emph{\textbf{v}})$.
\par
(c)\ If $u=\rho(\emph{\textbf{v}})$ then it is obvious that $G(t,u,\emph{\textbf{v}})\equiv\rho(\emph{\textbf{v}})$.
\par
Now we prove uniqueness of the solution. Since $B_{ D}(y,\emph{\textbf{v}})+\bar{B}_{ D}(y)$ satisfies Lipschitz condition with respect to $y$ for any fixed $(u,\emph{\textbf{v}})\in [0,1)\times[0,1)^D$, by ordinary differential equation theory, we know that (\ref{eq3.6}) has unique solution $y(t;u,\textbf{v})=G(t,u,\emph{\textbf{v}})$ for any fixed $(u,\emph{\textbf{v}})\in [0,1)\times[0,1)^D$.
\par
We now consider the case that $u=1$. Assume $\tilde{y}(t;1,\emph{\textbf{v}})$ is another solution of (\ref{eq3.6}). Since $\tilde{y}'(0;1,\emph{\textbf{v}})=B_{ D}(1)+\bar{B}_{ D}(1)<0$, we know that $\tilde{y}(t;1,\emph{\textbf{v}})\uparrow 1$ as $t\downarrow 0$. Hence, for any sufficient small $\varepsilon>0$, $\tilde{y}(\varepsilon;1,\emph{\textbf{v}})\in (\rho(\emph{\textbf{v}}),1)$. It is easy to see that $\bar{y}(t;1,\emph{\textbf{v}})=:\tilde{y}(\varepsilon+t;1,\emph{\textbf{v}})$ is a solution of (\ref{eq3.6}) with the initial condition $\bar{y}(0;1,\emph{\textbf{v}})=\tilde{y}(\varepsilon;1,\emph{\textbf{v}})$. Let
$$
\delta_{\varepsilon}=\inf\{t\geq 0; G(t,1,\emph{\textbf{v}})=\tilde{y}(\varepsilon;1,\emph{\textbf{v}})\}.
$$
Then, $\delta_{\varepsilon}\downarrow 0$ as $\varepsilon\downarrow 0$.
Similarly, $\hat{y}(t;1,\emph{\textbf{v}})=G(\delta_{\varepsilon}+t,1,\emph{\textbf{v}})$ is also a solution of (\ref{eq3.6}) with the initial condition $\hat{y}(0;1,\emph{\textbf{v}})=\tilde{y}(\varepsilon)$.
However, (\ref{eq3.6}) has unique solution with initial condition $\tilde{y}(\varepsilon)\in [0,1)$. Therefore,
 \begin{eqnarray*}
\tilde{y}(\varepsilon+t;1,\emph{\textbf{v}})=G(\delta_{\varepsilon}+t,1,\emph{\textbf{v}}), \ \ \forall t\geq 0,
\end{eqnarray*}
and hence, $\tilde{y}(t;1,\emph{\textbf{v}})\equiv G(t,1,\emph{\textbf{v}})$. The proof is complete. \hfill $\Box$
\end{proof}

\par
The following theorem gives the joint probability generating function of $( D-1)$-crossing numbers until time $t$, i.e., the joint probability generating function of $\emph{\textbf{Y}}(t)$.
\par
\begin{theorem}\label{th3.1}
Suppose that $\{X(t); t\geq 0\}$ is a Markov branching process with $X(0)=1$. Then
the joint probability generating function of $\textbf{Y}(t)$ is given by
\begin{eqnarray*}
E[\textbf{v}^{\textbf{Y}(t)}|X(0)=1]=G(t,1,\textbf{v}),\ \ \textbf{v}\in [0,1]^D,
\end{eqnarray*}
where $y=G(t,u,\textbf{v})$ is the unique solution of $(\ref{eq3.6})$. Furthermore,
\begin{eqnarray*}
P(\textbf{Y}(t)=\textbf{m}|X(0)=1)=g_{\textbf{m}}(t),\ \ \textbf{m}\in \mathbb{Z}_+^D,
\end{eqnarray*}
where
\begin{eqnarray*}
\begin{cases}
g_{\textbf{0}}(t)=G(t,1,\emph{\textbf{0}})\\
g_{\textbf{m}}(t)
=\bar{B}_{ D}(g_{\emph{\textbf{0}}}(t))\cdot
\int_0^t\frac{F_{\textbf{m}}(s)}{ \bar{B}_{ D}(g_{\textbf{0}}(s))}ds, \ \textbf{m}\neq \textbf{0}
\end{cases}
\end{eqnarray*}
with
\begin{eqnarray*}
F_{\textbf{m}}(t)=\sum_{i\in D}b_i\cdot g^{*(i)}_{\textbf{m}-\textbf{e}_i}(t)
+\sum_{i\in  D^c}b_i\cdot\sum_{\textbf{l}^{(1)},\cdots,\textbf{l}^{(i)}\neq \textbf{m},\  \textbf{l}^{(1)}+\cdots+\textbf{l}^{(i)}=\textbf{m}}g_{\textbf{l}^{(1)}}(t)\cdots g_{\textbf{l}^{(i)}}(t)
\end{eqnarray*}
and $\{g^{*(i)}_{\textbf{m}}(t);\ \textbf{m}\in \mathbb{Z}_+^D\}$ is the $i$'th convolution of $\{g_{\textbf{m}}(t);\ \textbf{m}\in \mathbb{Z}_+^D\}$.
\end{theorem}
\par
\begin{proof} Let $P(t)=(p_{(i,\emph{\textbf{m}}), (j,\emph{\textbf{l}})}(t); (i,\emph{\textbf{m}}), (j,\emph{\textbf{l}})\in \mathbb{Z}_+\times\mathbb{Z}_+^D)$ be the transition probability of $(X(t),\emph{\textbf{Y}}(t))$. We first prove that
\begin{eqnarray}\label{eq3.7}
G(t,u,\emph{\textbf{v}})=\sum_{(j,\emph{\textbf{l}})\in \mathbb{Z}_+\times\mathbb{Z}_+^D}p_{(1,\emph{\textbf{0}}),
 (j,\emph{\textbf{l}})}(t)u^j\emph{\textbf{v}}^{\emph{\textbf{l}}},\ \ (u,\emph{\textbf{v}})\in [0,1]\times [0,1)^D.
\end{eqnarray}
i.e., it suffices to prove that
$$
y(t,u,\emph{\textbf{v}})=\sum_{(j,\emph{\textbf{l}})\in \mathbb{Z}_+\times\mathbb{Z}_+^D}p_{(1,\emph{\textbf{0}}),
 (j,\emph{\textbf{l}})}(t)u^j\emph{\textbf{v}}^{\emph{\textbf{l}}}
$$
is the solution of (\ref{eq3.6}). Indeed, it follows from Kolmogorov backward equation that
\begin{eqnarray*}
&& p'_{(1,\emph{\textbf{0}}),
 (j,\emph{\textbf{l}})}(t)\\
 &=&\sum_{k\geq 0,\emph{\textbf{r}}\in \mathbb{Z}_+^D}q_{(1,\emph{\textbf{0}}),
 (k,\emph{\textbf{r}})}\cdot p_{(k,\emph{\textbf{r}}),(j,\emph{\textbf{l}})}(t)\\
 &=&\sum_{k\in  D}b_{k}\cdot p_{(k,\emph{\textbf{e}}_{k}),(j,\emph{\textbf{l}})}(t)+\sum_{k\in D^c}b_{k}\cdot p_{(k,\emph{\textbf{0}}),(j,\emph{\textbf{l}})}(t),\ \ \forall t\geq 0.
\end{eqnarray*}
Multiplying $u^j\emph{\textbf{v}}^{\emph{\textbf{l}}}$ on both sides of the above equality, then summing over $j$ and $\emph{\textbf{l}}$ and using Lemma~\ref{le3.1} yield that
\par
\begin{eqnarray*}
 && \sum_{(j,\emph{\textbf{l}})\in \mathbb{Z}_+\times\mathbb{Z}_+^D}p'_{(1,\textbf{0}),
 (j,\emph{\textbf{l}})}(t)\cdot u^j\emph{\textbf{v}}^{\emph{\textbf{l}}}\\
 &=&\sum_{k\in  D}b_{k}\cdot \sum_{(j,\emph{\textbf{l}})\in \mathbb{Z}_+\times\mathbb{Z}_+^D}p_{(k,\emph{\textbf{e}}_{k}),
 (j,\emph{\textbf{l}})}(t)\cdot u^j\emph{\textbf{v}}^{\emph{\textbf{l}}}+\sum_{k\in  D^c}b_{k}\cdot \sum_{(j,\emph{\textbf{l}})\in \mathbb{Z}_+\times\mathbb{Z}_+^D}p_{(k,\emph{\textbf{0}}),
 (j,\emph{\textbf{l}})}(t)\cdot u^j\emph{\textbf{v}}^{\emph{\textbf{l}}}\\
 &=&\sum_{k\in  D}b_{k}\cdot [\sum_{(j,\emph{\textbf{l}})\in \mathbb{Z}_+\times\mathbb{Z}_+^D}p_{(1,\emph{\textbf{0}}),
 (j,\emph{\textbf{l}})}(t)\cdot u^j\emph{\textbf{v}}^{\emph{\textbf{l}}}]^k\cdot\emph{\textbf{v}}
 ^{\emph{\textbf{e}}_k}+\sum_{k\in  D^c}b_{k}\cdot [\sum_{(j,\emph{\textbf{l}})\in \mathbb{Z}_+\times\mathbb{Z}_+^D}p_{(1,\emph{\textbf{0}}),
 (j,\emph{\textbf{l}})}(t)\cdot u^j\emph{\textbf{v}}^{\emph{\textbf{l}}}]^k\\
 &=& B_{ D}(\sum_{(j,\emph{\textbf{l}})\in \mathbb{Z}_+\times\mathbb{Z}_+^D}p_{(1,\emph{\textbf{0}}),
 (j,\emph{\textbf{l}})}(t)\cdot u^j\emph{\textbf{v}}^{\emph{\textbf{l}}},\emph{\textbf{v}})+\bar{B}_{ D}
 (\sum_{(j,\emph{\textbf{l}})\in \mathbb{Z}_+\times\mathbb{Z}_+^D}p_{(1,\emph{\textbf{0}}),
 (j,\emph{\textbf{l}})}(t)\cdot u^j\emph{\textbf{v}}^{\emph{\textbf{l}}})
\end{eqnarray*}
which implies that $y(t,u,\emph{\textbf{v}})=\sum_{(j,\emph{\textbf{l}})\in \mathbb{Z}_+\times\mathbb{Z}_+^D}p_{(1,\emph{\textbf{0}}),
 (j,\emph{\textbf{l}})}(t)u^j\emph{\textbf{v}}^{\emph{\textbf{l}}}$ satisfies
$$
\frac{\partial y}{\partial t}=B_{ D}(y,\emph{\textbf{v}})+\bar{B}_{ D}(y).
$$
Finally, it is easy to see that
$$
y(0,u,\emph{\textbf{v}})=u.
$$
Therefore (\ref{eq3.7}) is proved. Hence, it follows from (\ref{eq3.7}) and $\emph{\textbf{Y}}(0)=\emph{\textbf{0}}$ that
\begin{eqnarray*}
E[\emph{\textbf{v}}^{\emph{\textbf{Y}}(t)}|X(0)=1]
&=& \sum_{\emph{\textbf{l}}\in\mathbb{Z}_+^D}P(\emph{\textbf{Y}}(t)=\emph{\textbf{l}}
|X(0)=1,\emph{\textbf{Y}}(0)=\emph{\textbf{0}})
\cdot\emph{\textbf{v}}^{\emph{\textbf{l}}}\\
 &=& \sum_{\emph{\textbf{l}}\in \mathbb{Z}_+^D}\sum_{j=0}^{\infty}p_{(1,\emph{\textbf{0}}),
 (j,\emph{\textbf{l}})}(t)\cdot\emph{\textbf{v}}^{\emph{\textbf{l}}}\\
 &=& G(t,1,\emph{\textbf{v}}).
\end{eqnarray*}
\par
Finally, it follows from the above proof that $G(t,u,\emph{\textbf{v}})$ can be expanded as a multivariate nonnegative Taylor series. Suppose that
\begin{eqnarray*}
  G(t,1,\emph{\textbf{v}})=\sum_{\emph{\textbf{l}}\in \mathbb{Z}_+^D}g_{\emph{\textbf{l}}}(t)\emph{\textbf{v}}^{\emph{\textbf{l}}}.
\end{eqnarray*}
By (\ref{eq3.6}),
\begin{eqnarray*}
  \sum_{\emph{\textbf{l}}\in \mathbb{Z}_+^D}g'_{\emph{\textbf{l}}}(t)\emph{\textbf{v}}^{\emph{\textbf{l}}}
  &=&\sum_{i\in D}b_i\cdot[\sum_{\emph{\textbf{l}}\in \mathbb{Z}_+^D}g_{\emph{\textbf{l}}}(t)\emph{\textbf{v}}^
  {\emph{\textbf{l}}}
  ]^i\cdot
  \emph{\textbf{v}}^{\emph{\textbf{e}}_i}+\sum_{i\in D^c}b_i\cdot[\sum_{\emph{\textbf{l}}\in \mathbb{Z}_+^D}g_{\emph{\textbf{l}}}(t)\emph{\textbf{v}}^{\emph{\textbf{l}}}
 ]^i\\
&=& \sum_{i\in D}b_i\sum_{\emph{\textbf{l}}\in \mathbb{Z}_+^D}g^{*(i)}_{\emph{\textbf{l}}}(t)\emph{\textbf{v}}^{\emph{\textbf{l}}
+\emph{\textbf{e}}_i}
+\sum_{i\in  D^c}b_i\sum_{\emph{\textbf{l}}\in \mathbb{Z}_+^D}g^{*(i)}_{\emph{\textbf{l}}}(t)\emph{\textbf{v}}
^{\emph{\textbf{l}}}\\
&=& \sum_{\emph{\textbf{l}}\in \mathbb{Z}_+^D}\sum_{i\in D}b_ig^{*(i)}_{\emph{\textbf{l}}}(t)
\emph{\textbf{v}}^{\emph{\textbf{l}}+\emph{\textbf{e}}_i}
+\sum_{\emph{\textbf{l}}\in \mathbb{Z}_+^D}\sum_{i\in D^c}b_jg^{*(i)}_{\emph{\textbf{l}}}(t)\emph{\textbf{v}}^{\emph{\textbf{l}}}\\
&=& \sum_{\emph{\textbf{l}}\in \mathbb{Z}_+^D\setminus \{\emph{\textbf{0}}\}}\sum_{i\in D}b_ig^{*(i)}_{\emph{\textbf{l}}
-\emph{\textbf{e}}_i}(t)
\emph{\textbf{v}}^{\emph{\textbf{l}}}
+\sum_{\emph{\textbf{l}}\in \mathbb{Z}_+^D}\sum_{i\in D^c}b_ig^{*(i)}_{\emph{\textbf{l}}}(t)\emph{\textbf{v}}^{\emph{\textbf{l}}}.
\end{eqnarray*}
where $\{g^{*(i)}_{\emph{\textbf{l}}}(t);\ \emph{\textbf{l}}\in \mathbb{Z}_+^D\}$ is the $i$th convolution of $\{g_{\emph{\textbf{l}}}(t);\ \emph{\textbf{l}}\in \mathbb{Z}_+^D\}$ and here we have used the notation $g^{*(i)}_{\emph{\textbf{l}}}(t)=0$ if $\emph{\textbf{l}}\notin \mathbb{Z}_+^D$.
Comparing the coefficients on the both sides of the above equality yields that
\begin{eqnarray}\label{eq3.8}
\begin{cases}
g'_{\emph{\textbf{0}}}(t)=\sum_{i\in D^c}b_ig^i_{\emph{\textbf{0}}}(t)=\bar{B}_{ D}
(g_{\emph{\textbf{0}}}(t)),\\
g'_{\emph{\textbf{l}}}(t)=\sum_{i\in D}b_ig^{*(i)}_{\emph{\textbf{l}}
-\emph{\textbf{e}}_i}(t)
+\sum_{i\in  D^c}b_ig^{*(i)}_{\emph{\textbf{l}}}(t),\ \ \emph{\textbf{l}}\neq \emph{\textbf{0}}.
\end{cases}
\end{eqnarray}
It is easy to see that $g_{\emph{\textbf{l}}}(0)=P(\emph{\textbf{Y}}(0)=\emph{\textbf{l}}|X(0)=1)
=\delta_{\emph{\textbf{0}},\emph{\textbf{l}}}$ and hence $g_{\emph{\textbf{0}}}(t)=G(t,1,\emph{\textbf{0}})$. On the other hand, by the second equation of (\ref{eq3.8}),
\begin{eqnarray*}
g'_{\emph{\textbf{l}}}(t)-g_{\emph{\textbf{l}}}(t)\bar{B}'_{ D}
(g_{\emph{\textbf{0}}}(t))=F_{\emph{\textbf{l}}}(t).
\end{eqnarray*}
Therefore, note that $g_{\emph{\textbf{k}}}(0)=P(\emph{\textbf{Y}}(0)=\emph{\textbf{k}}|X(0)=1)=0$ for all $\emph{\textbf{k}}\neq \emph{\textbf{0}}$, we have
\begin{eqnarray*}
g_{\emph{\textbf{l}}}(t)e^{-\int_0^t\bar{B}'_{ D}
(g_{\emph{\textbf{0}}}(s))ds}
=\int_0^tF_{\emph{\textbf{l}}}(s)\cdot e^{-\int_0^s\bar{B}'_{ D}(g_{\emph{\textbf{0}}}(x))dx}ds.
\end{eqnarray*}
Note that
\begin{eqnarray*}
\int_0^s\bar{B}'_{ D}(g_{\emph{\textbf{0}}}(x))dx=\int_0^s\bar{B}'_{ D}(g_{\emph{\textbf{0}}}(x))\cdot \frac{g'_{\emph{\textbf{0}}}(x)}{\bar{B}_{ D}
(g_{\emph{\textbf{0}}}(x))}dx=\ln \bar{B}_{ D}
(g_{\emph{\textbf{0}}}(s))
\end{eqnarray*}
Hence,
\begin{eqnarray*}
g_{\emph{\textbf{l}}}(t)
=\bar{B}_{ D}(g_{\emph{\textbf{0}}}(t))\cdot
\int_0^t\frac{F_{\emph{\textbf{l}}}(s)}{ \bar{B}_{ D}(g_{\emph{\textbf{0}}}(s))}ds,\ \ \emph{\textbf{l}}\neq \emph{\textbf{0}}.
\end{eqnarray*}
The proof is complete. \hfill $\Box$
\end{proof}
\par
\begin{remark}\label{re3.1}
(i)\ Generally, if $X(t)$ starts from $X(0)=i (>1)$, then the joint probability generating function of $( D-1)$-crossing numbers until time $t$ is
\begin{eqnarray*}
E[\emph{\textbf{v}}^{\emph{\textbf{Y}}(t)}|X(0)=i]=[G(t,1,\emph{\textbf{v}})]^i.
\end{eqnarray*}
\par
(ii)\ By carefully checking the proof of Theorem~\ref{th3.1}, we see that $G(t,u,\emph{\textbf{v}})$ can be expanded as a nonnegative multivariate Taylor series
 \begin{eqnarray*}
G(t,u,\emph{\textbf{v}})=\sum_{(j,\emph{\textbf{l}})\in \mathbb{Z}_+\times\mathbb{Z}_+^D}g_{
 j,\emph{\textbf{l}}}(t)u^j\emph{\textbf{v}}^{\emph{\textbf{l}}},\ \ (u,\emph{\textbf{v}})\in [0,1]\times [0,1)^D,
\end{eqnarray*}
where $g_{j,\emph{\textbf{l}}}(t)=p_{(1,\emph{\textbf{0}}),(j,\emph{\textbf{l}})}(t)$ for any $(j,\emph{\textbf{l}})\in \mathbb{Z}_+\times\mathbb{Z}_+^D$. Therefore, if the solution $G(t,u,\emph{\textbf{v}})$ is known, then we can obtain $\{p_{(1,\emph{\textbf{0}}),
 (j,\emph{\textbf{l}})}(t); (j,\emph{\textbf{l}})\in \mathbb{Z}_+\times\mathbb{Z}_+^D\}$. Hence,  by (\ref{eq3.5}), the transition probability function $P(t)=(p_{(i,\emph{\textbf{m}}), (j,\emph{\textbf{l}})}(t); (i,\emph{\textbf{m}}), (j,\emph{\textbf{l}})\in \mathbb{Z}_+\times\mathbb{Z}_+^D)$ can be obtained.
\end{remark}
\par
The following theorem gives a recursive algorithm of $g_{j,\emph{\textbf{l}}}(t)$.
\par
\begin{theorem}\label{th3.2}
Suppose that $\{X(t); t\geq 0\}$ is a Markov branching process with $X(0)=1$.
\par
\rm{(i)}\ If $0\in  D$ or $b_0=0$, then
$g_{j\emph{\textbf{l}}}(t)$ is given by
\begin{eqnarray*}
\begin{cases}
g_{0\emph{\textbf{0}}}(t)=0\\
g_{j\emph{\textbf{l}}}(t)=e^{b_1t}[\delta_{j,1}\delta_{\emph{\textbf{l}},
\emph{\textbf{0}}}+\int_0^t F_{j,\emph{\textbf{l}}}(s)e^{-b_1s}ds],\ \ (j,\emph{\textbf{l}})\neq (0,\emph{\textbf{0}}),
\end{cases}
\end{eqnarray*}
\par
\rm{(ii)}\ If $0\notin  D$ and $b_0>0$, then
$g_{j\emph{\textbf{l}}}(t)$ is given by
\begin{eqnarray*}
\begin{cases}
g_{0\emph{\textbf{0}}}(t)=G(t,0,\emph{\textbf{0}})\\
g_{j,\emph{\textbf{l}}}(t)=\bar{B}_{ D}(g_{0 \emph{\textbf{0}}}(t))\cdot[\delta_{j,1}\delta_{\emph{\textbf{l}},\emph{\textbf{0}}}
b_0^{-1}
+\int_0^t\frac{F_{j,\emph{\textbf{l}}}(s)}{\bar{B}_{ D}(g_{0 \emph{\textbf{0}}}(s))}ds],\ \ (j,\emph{\textbf{l}})\neq (0,\emph{\textbf{0}}),
\end{cases}
\end{eqnarray*}
where
\begin{eqnarray*}
F_{j,\emph{\textbf{l}}}(t)=\sum_{i\in  D}b_i\cdot g^{*(i)}_{j \emph{\textbf{l}}-\emph{\textbf{e}}_i}(t)+\sum_{i\in D^c}b_i\cdot \sum_{(j_1,\emph{\textbf{l}}^{(1)}),\cdots,
(j_i,\emph{\textbf{l}}^{(i)})\neq (j,\emph{\textbf{l}}),\sum_{k=1}^i(j_k,\emph{\textbf{l}}^{(k)})
=(j,\emph{\textbf{l}})}
g_{j_1 \emph{\textbf{l}}^{(1)}}(t)\cdots g_{j_i \emph{\textbf{l}}^{(i)}}(t).
\end{eqnarray*}
and $\{g^{*(i)}_{j \emph{\textbf{l}}}(t);\ (j,\emph{\textbf{l}})\in \mathbb{Z}_+\times\mathbb{Z}_+^D\}$ is the $i$'th convolution of $\{g_{j \emph{\textbf{l}}}(t);\ (j,\emph{\textbf{l}})\in \mathbb{Z}_+\times\mathbb{Z}_+^D\}$. Here $g^{*(i)}_{j \emph{\textbf{l}}}(t)=0$ if $\emph{\textbf{l}}\notin \mathbb{Z}_+^D$.
\end{theorem}
\par
\begin{proof}
By Remark~\ref{re3.1}, $G(t,u,\emph{\textbf{v}})$ can be expanded as a nonnegative multivariate Taylor series
\begin{eqnarray*}
  G(t,u,\emph{\textbf{v}})=\sum_{(j,\emph{\textbf{l}})\in \mathbb{Z}_+\times\mathbb{Z}_+^D}g_{j \emph{\textbf{l}}}(t)u^j\emph{\textbf{v}}^{\emph{\textbf{l}}}.
\end{eqnarray*}
In order to get the coefficients $g_{j \emph{\textbf{l}}}(t)$, substitute $G(t,u,\emph{\textbf{v}})$ into (\ref{eq3.6}), we obtain
\begin{eqnarray*}
  &&\sum_{(j,\emph{\textbf{l}})\in \mathbb{Z}_+\times\mathbb{Z}_+^D}g'_{j \emph{\textbf{l}}}(t)u^j\emph{\textbf{v}}^{\emph{\textbf{l}}}\\
  &=&\sum_{i\in D}b_i\cdot[\sum_{(j,\emph{\textbf{l}})\in \mathbb{Z}_+\times\mathbb{Z}_+^D}g_{j \emph{\textbf{l}}}(t)u^j\emph{\textbf{v}}^{\emph{\textbf{l}}}]^i
  \cdot\emph{\textbf{v}}^{\emph{\textbf{e}}_i}+\sum_{i\in D^c}b_i\cdot[\sum_{(j,\emph{\textbf{l}})\in \mathbb{Z}_+\times\mathbb{Z}_+^D}g_{j \emph{\textbf{l}}}(t)u^j\emph{\textbf{v}}^{\emph{\textbf{l}}}]^i\\
&=& \sum_{i\in D}b_i\sum_{(j,\emph{\textbf{l}})\in \mathbb{Z}_+\times\mathbb{Z}_+^D}g^{*(i)}_{j \emph{\textbf{l}}}(t)u^j\emph{\textbf{v}}^{\emph{\textbf{l}}
+\emph{\textbf{e}}_i}
+\sum_{i\in  D^c}b_i\sum_{(j,\emph{\textbf{l}})\in \mathbb{Z}_+\times\mathbb{Z}_+^D}g^{*(i)}_{j \emph{\textbf{l}}}(t)u^j\emph{\textbf{v}}^{\emph{\textbf{l}}}\\
&=& \sum_{(j,\emph{\textbf{l}})\in \mathbb{Z}_+\times\mathbb{Z}_+^D}\sum_{i\in D}b_ig^{*(i)}_{j \emph{\textbf{l}}}(t)u^j\emph{\textbf{v}}^{\emph{\textbf{l}}
+\emph{\textbf{e}}_i}
+\sum_{(j,\emph{\textbf{l}})\in \mathbb{Z}_+\times\mathbb{Z}_+^D}\sum_{i\in  D^c}b_ig^{*(i)}_{j \emph{\textbf{l}}}(t)u^j\emph{\textbf{v}}^{\emph{\textbf{l}}}\\
&=& \sum_{(j,\emph{\textbf{l}})\in \mathbb{Z}_+\times\mathbb{Z}_+^D\setminus \{\emph{\textbf{0}}\}}\sum_{i\in D}b_ig^{*(i)}_{j \emph{\textbf{l}}-\emph{\textbf{e}}_i}(t)
u^j\emph{\textbf{v}}^{\emph{\textbf{l}}}
+\sum_{(j,\emph{\textbf{l}})\in \mathbb{Z}_+\times\mathbb{Z}_+^D}\sum_{i\in  D^c}b_ig^{*(i)}_{j \emph{\textbf{l}}}(t)u^j\emph{\textbf{v}}^{\emph{\textbf{l}}},
\end{eqnarray*}
 here we have used the notation $g^{*(i)}_{j \emph{\textbf{l}}}(t)=0$ if $\emph{\textbf{l}}\notin \mathbb{Z}_+^D$.
Comparing the coefficients on the both sides of the above equality yields that
\begin{eqnarray}\label{eq3.10}
g'_{j \emph{\textbf{l}}}(t)=\sum_{i\in D}b_ig^{*(i)}_{j \emph{\textbf{l}}-\emph{\textbf{e}}_i}(t)
+\sum_{i\in  D^c}b_ig^{*(i)}_{j \emph{\textbf{l}}}(t),\ \ (j,\emph{\textbf{l}})\in \mathbb{Z}_+\times\mathbb{Z}_+^D.
\end{eqnarray}
\par
It is easy to see that
$$
g_{j \emph{\textbf{0}}}(0)=P(X(0)=j, \emph{\textbf{Y}}(0)=0|X(0)=1)=\delta_{1j}.
$$
\par
For $(j,\emph{\textbf{l}})=(0,\emph{\textbf{0}})$, by (\ref{eq3.10}),
 \begin{eqnarray*}
g'_{0 \emph{\textbf{0}}}(t)=\sum_{i\in D^c}b_ig^i_{0,\emph{\textbf{0}}}(t)=\bar{B}_{ D}
(g_{0,\emph{\textbf{0}}}(t)),
\end{eqnarray*}
which implies
 $$
 g_{0 \emph{\textbf{0}}}(t)=G(t,0,\emph{\textbf{0}}).
$$
For $(j,\emph{\textbf{l}})\neq(0,\emph{\textbf{0}})$, by (\ref{eq3.10}),
\begin{eqnarray}\label{eq3.11}
&&g'_{j \emph{\textbf{l}}}(t)\nonumber\\
&=&\sum_{i\in  D}b_ig^{*(i)}_{j \emph{\textbf{l}}-\emph{\textbf{e}}_i}(t)
+\sum_{i\in  D^c}b_ig^{*(i)}_{j \emph{\textbf{l}}}(t)\nonumber\\
&=&g_{j \emph{\textbf{l}}}(t)\sum_{i\in  D^c}ib_ig^{i-1}_{0 \emph{\textbf{0}}}(t)\nonumber\\
&&+\sum_{i\in D}b_ig^{*(i)}_{j \emph{\textbf{l}}-\emph{\textbf{e}}_i}(t)
+\sum_{i\in  D^c}b_i\sum_{(j_1,\emph{\textbf{l}}^{(1)}),\cdots,(j_i,
\emph{\textbf{l}}^{(i)})
\neq (j,\emph{\textbf{l}}),\sum_{k=1}^i(j_k,\emph{\textbf{l}}^{(k)})=
(j,\emph{\textbf{l}})}
g_{j_1 \emph{\textbf{l}}^{(1)}}(t)\cdots g_{j_i \emph{\textbf{l}}^{(i)}}(t)\nonumber\\
&=&g_{j \emph{\textbf{l}}}(t)\bar{B}'_{ D}(g_{0 \emph{\textbf{0}}}(t))+F_{j,\emph{\textbf{l}}}(t).
\end{eqnarray}
\par
If $0\in D$ or $b_0=0$, then $g_{0 \emph{\textbf{0}}}(t)=p_{(1,\emph{\textbf{0}}),(0,\emph{\textbf{0}})}(t)=0$
and hence
$$
\bar{B}'_{ D}(g_{0 \emph{\textbf{0}}}(t))=b_1.
$$
By (\ref{eq3.11}),
\begin{eqnarray*}
&&g'_{j \emph{\textbf{l}}}(t)=b_1g_{j \emph{\textbf{l}}}(t)+F_{j,\emph{\textbf{l}}}(t).
\end{eqnarray*}
Hence,
\begin{eqnarray*}
g_{j \emph{\textbf{l}}}(t)=e^{b_1t}[\delta_{j,1}\delta_{\emph{\textbf{l}},
\emph{\textbf{0}}}+\int_0^t F_{j,\emph{\textbf{l}}}(s)e^{-b_1s}ds].
\end{eqnarray*}
\par
If $0\notin D$ and $b_0>0$, then
$$
e^{\int_0^t\bar{B}'_{ D}(g_{0 \emph{\textbf{0}}}(s))ds}
=e^{\int_0^t\bar{B}'_{ D}(g_{0 \emph{\textbf{0}}}(s))\cdot\frac{g'_{0\emph{\textbf{0}}}(s)}{\bar{B}_{ D}
(g_{0 \emph{\textbf{0}}}(s))}ds}
=\frac{\bar{B}_{ D}(g_{0\emph{\textbf{0}}}(t))}
{b_0}
$$
Hence,
\begin{eqnarray*}
g_{j \emph{\textbf{l}}}(t)=
\bar{B}_{ D}(g_{0 \emph{\textbf{0}}}(t))\cdot[\delta_{j,1}\delta_{\emph{\textbf{l}},\emph{\textbf{0}}}
b_0^{-1}
+\int_0^t\frac{F_{j,\emph{\textbf{l}}}(s)}{\bar{B}_{ D}(g_{0 \emph{\textbf{0}}}(s))}ds],\ \ (j,\emph{\textbf{l}})\neq(0,\emph{\textbf{0}}).
\end{eqnarray*}

The proof is complete. \hfill $\Box$
\end{proof}
\par
As direct consequences of Theorem~\ref{th3.1} and Remark~\ref{re3.1}, the following corollaries give the probability distributions of death number and $(m-1)$-range up crossing number until time $t$ for fixed $m>1$.
\par
\begin{corollary}\label{cor3.1}
Suppose that $\{X(t); t\geq 0\}$ is a Markov branching process with $X(0)=1$. Then
the probability generating function of death number until time $t$ is given by
\begin{eqnarray*}
E[v^{Y_0(t)}|X(0)=1]=G(t,1,v),\ \ v\in [0,1],
\end{eqnarray*}
where $G(t,u,v)$ is the unique solution of the equation
\begin{eqnarray*}
\begin{cases}\frac{\partial y}{\partial t}=B(y)-b_0(1-v),\\
y|_{t=0}=u,
\end{cases} u,v\in [0,1].
\end{eqnarray*}
\end{corollary}

\par
\begin{corollary}\label{cor3.2}
Suppose that $\{X(t); t\geq 0\}$ is a Markov branching process with $X(0)=1$ and $m(>1)$ is fixed. Then
the probability generating function of $(m-1)$-range up-crossing number until time $t$ is given by
\begin{eqnarray*}
E[v^{Y_m(t)}|X(0)=1]=G(t,1,v),\ \ v\in [0,1],
\end{eqnarray*}
where $G(t,u,v)$ is the unique solution of the equation
\begin{eqnarray*}
\begin{cases}\frac{\partial y}{\partial t}=B(y)-b_m(1-v)y^m,\\
y|_{t=0}=u,
\end{cases} u,v\in [0,1].
\end{eqnarray*}
\end{corollary}
\par
Now, we consider the jumping property of $\{X(t);t\geq 0\}$ until extinction. Let
\begin{eqnarray*}
\tau=\inf\{t\geq 0;\ X(t)=0\}
\end{eqnarray*}
be the extinction time of $X(t)$.
\par
By Theorem~\ref{th3.1}, we can get the following result which is due to Li Y. and Li J.~\cite{LLC20}.
\par
\begin{theorem}\label{th3.3}
Suppose that $\{X(t); t\geq 0\}$ is a Markov branching process with $X(0)=1$. Then
the probability generating function $G(\textbf{v})$ of $( D-1)$-range crossing numbers conditioned on $\tau<\infty$ is given by
\begin{eqnarray*}
G(\textbf{v})=\rho^{-1}\cdot G(\infty,1,\textbf{v})=\rho^{-1}\cdot \rho(\textbf{v}),\ \ \textbf{v}\in [0,1]^D,
\end{eqnarray*}
where $\rho$ is the minimal nonnegative root of $B(u)=0$. Furthermore, if $\rho<1$ then
\begin{eqnarray}\label{eq3.12}
P(\textbf{Y}(\infty)={\bf \infty}|\tau=\infty)=1.
\end{eqnarray}
\end{theorem}
\par
\begin{proof}
It follows from Theorem~\ref{th3.1} that
\begin{eqnarray}\label{eq3.13}
G(t,1,\emph{\textbf{v}})=\sum_{\emph{\textbf{l}}\in \mathbb{Z}_+^D}p_{(1,\emph{\textbf{0}}),(0,\emph{\textbf{l}})}(t)
\emph{\textbf{v}}^{\emph{\textbf{l}}}+\sum_{\emph{\textbf{l}}\in \mathbb{Z}_+^D}(\sum_{j=1}^{\infty}p_{(1,\emph{\textbf{0}}),(j,\emph{\textbf{l}})}
(t))
\emph{\textbf{v}}^{\emph{\textbf{l}}},\ \ \ \forall t\geq 0.
\end{eqnarray}
By (\ref{eq3.3}) with $i=1$ and $u=\rho(\emph{\textbf{v}})$,
\begin{eqnarray}\label{eq3.14}
\rho(\emph{\textbf{v}})=\sum_{\emph{\textbf{l}}\in \mathbb{Z}_+^D}p_{(1,\emph{\textbf{0}}),(0,\emph{\textbf{l}})}(t)
\emph{\textbf{v}}^{\emph{\textbf{l}}}+\sum_{\emph{\textbf{l}}\in \mathbb{Z}_+^D}(\sum_{j=1}^{\infty}
p_{(1,\emph{\textbf{0}}),(j,\emph{\textbf{l}})}(t)
\rho(\emph{\textbf{v}})^j)
\emph{\textbf{v}}^{\emph{\textbf{l}}},\ \ \forall t\geq 0.
\end{eqnarray}
Noting that $(0,\emph{\textbf{l}})$ is absorbing state for any $\emph{\textbf{l}}$, we can let $t\rightarrow \infty$ in (\ref{eq3.13}) and (\ref{eq3.14}) to get
\begin{eqnarray}\label{eq3.15}
G(\infty,1,\emph{\textbf{v}})=\sum_{\emph{\textbf{l}}\in \mathbb{Z}_+^D}p_{(1,\emph{\textbf{0}}),(0,\emph{\textbf{l}})}(\infty)
\emph{\textbf{v}}^{\emph{\textbf{l}}}+\lim_{t\rightarrow \infty}\sum_{\emph{\textbf{l}}\in \mathbb{Z}_+^D}(\sum_{j=1}^{\infty}p_{(1,\emph{\textbf{0}}),(j,\emph{\textbf{l}})}
(t))
\emph{\textbf{v}}^{\emph{\textbf{l}}}
\end{eqnarray}
and
\begin{eqnarray}\label{eq3.16}
\rho(\emph{\textbf{v}})=\sum_{\emph{\textbf{l}}\in \mathbb{Z}_+^D}p_{(1,\emph{\textbf{0}}),(0,\emph{\textbf{l}})}(\infty)
\emph{\textbf{v}}^{\emph{\textbf{l}}}.
\end{eqnarray}
By (\ref{eq3.15}) and (\ref{eq3.16}),
\begin{eqnarray*}
\sum_{\emph{\textbf{l}}\in \mathbb{Z}_+^D}p_{(1,\emph{\textbf{0}}),(0,\emph{\textbf{l}})}(\infty)
\emph{\textbf{v}}^{\emph{\textbf{l}}}=G(\infty,1,\emph{\textbf{v}})=
\rho(\emph{\textbf{v}}),
\end{eqnarray*}
Therefore,
\begin{eqnarray*}
G(\emph{\textbf{v}})&=&\sum_{\emph{\textbf{l}}\in \mathbb{Z}_+^D}
P(\emph{\textbf{Y}}(\tau)=\emph{\textbf{l}}\ |\tau<\infty)\cdot \emph{\textbf{v}}^{\emph{\textbf{l}}}\\
&=&\rho^{-1}\cdot\sum_{\emph{\textbf{l}}\in \mathbb{Z}_+^D}p_{(1,\emph{\textbf{0}}),(0,\emph{\textbf{l}})}(\infty)
\emph{\textbf{v}}^{\emph{\textbf{l}}}\\
&=& \rho^{-1}\cdot G(\infty,1,\emph{\textbf{v}})\\
&=& \rho^{-1}\cdot\rho(\emph{\textbf{v}}).
\end{eqnarray*}
Again by (\ref{eq3.15}) and (\ref{eq3.16}),
\begin{eqnarray*}
&&P(\emph{\textbf{Y}}(\tau)=\emph{\textbf{l}}\ |\tau=\infty)\\
&\leq & P(\emph{\textbf{Y}}(\tau)\leq \emph{\textbf{l}}\ |\tau=\infty)\\
&=& (1-\rho)^{-1}\cdot P(\emph{\textbf{Y}}(\tau)\leq \emph{\textbf{l}},\tau=\infty)\\
&=& (1-\rho)^{-1}\cdot \lim_{t\rightarrow \infty}P(\emph{\textbf{Y}}(t)\leq \emph{\textbf{l}},\tau>t)\\
&=& (1-\rho)^{-1}\cdot \lim_{t\rightarrow \infty}\sum_{\emph{\textbf{m}}\leq \emph{\textbf{l}}}\sum_{j=1}^{\infty}p_{(1,\emph{\textbf{0}}),
(j,\emph{\textbf{m}})}(t)\\
&=&0,
\end{eqnarray*}
which implies (\ref{eq3.12}), where $\emph{\textbf{m}}\leq \emph{\textbf{l}}$ means every component of $\emph{\textbf{m}}$ is not bigger than $\emph{\textbf{l}}$.
The proof is complete.
\hfill $\Box$
\end{proof}
\par
It can be seen from Theorems~\ref{th3.1} and \ref{th3.3}, in order to obtain the joint probability generating function of $\emph{\textbf{Y}}(t)$, the key role is to find the solution of (\ref{eq3.6}). Therefore, we consider two special and important cases to illustrate how to find the solution of (\ref{eq3.6}) and Taylor expansion of the solution.
\par
\begin{corollary}\label{cor3.3}
Suppose that $X(t)$ is a birth-death type branching process with death rate $pb$ and birth rate $qb$ $($where $b>0,\ p\in(0,1), p+q=1$$)$, $X(0)=1$. $Y(t)$ is the death number of $\{X(t);t\geq 0\}$ until time $t$. Then
\par
{\rm{(i)}}\ the probability generating function of $Y(t)$ is given by
\begin{eqnarray*}
E[v^{Y(t)}]=G(t,1,v)=\beta(v)+\frac{\alpha(v)-\beta(v)}{1+\frac{\alpha(v)-1}{1-\beta(v)}\cdot e^{(\alpha(v)-\beta(v))bqt}},
\end{eqnarray*}
where
\begin{eqnarray*}
\alpha(v)=\frac{1+\sqrt{1-4pqv}}{2q},\ \ \ \ \beta(v)&=&\frac{1-\sqrt{1-4pqv}}{2q}.
\end{eqnarray*}
More specifically, $g_n(t)=P(Y(t)=n)\ (n\geq 0)$ is given by
\begin{eqnarray*}
\begin{cases}
g_0(t)=\frac{1}{q+p\cdot e^{bt}},\\
g_n(t)
=\frac{e^{bt}}{(q+p\cdot e^{bt})^2}\cdot
\int_0^t(q+p\cdot e^{bs})^2e^{-bs}F_n(s)ds, \ n\geq 1
\end{cases}
\end{eqnarray*}
with
$$
F_n(t)=bp\delta_{1,n}+bq\sum_{k=1}^{n-1}g_k(t)g_{n-k}(t).
$$
\par
{\rm{(ii)}}\ The probability generating function of death number conditioned on $\tau<\infty$ is given by
\begin{eqnarray*}
E[v^{Y(\tau)}|\tau<\infty]=\beta(v),
\end{eqnarray*}
where $\tau=\inf\{t\geq 0; X(t)=0\}$ is the extinction time of $\{X(t);t\geq 0\}$.
\end{corollary}
\begin{proof} Note that $ D=\{0\}$, $B(y)=b(p-y+qy^2)$ and
$$
B_{ D}(y,v)+\bar{B}_{ D}(y)=b(pv-y+qy^2).
$$
It is obvious that for any $v\in [0,1)$,
$$
B_{ D}(y,v)+\bar{B}_{ D}(y)=b(pv-y+qy^2)=bq(y-\alpha(v))(y-\beta(v)),
$$
where
\begin{eqnarray*}
\alpha(v)=\frac{1+\sqrt{1-4pqv}}{2q},\ \ \ \ \beta(v)&=&\frac{1-\sqrt{1-4pqv}}{2q}.
\end{eqnarray*}
Therefore, (\ref{eq3.6}) becomes
\begin{eqnarray}\label{eq3.17}
  \begin{cases}\frac{dy}{(y-\alpha(v))(y-\beta(v))}=bqdt,\\
  y(0)=u.
  \end{cases}
\end{eqnarray}
Note that $\alpha(v)>\beta(v)$ for $v\in [0,1)$. Solve (\ref{eq3.17}), one get
\begin{eqnarray*}
y(t)=G(t,u,v)=\beta(v)+\frac{\alpha(v)-\beta(v)}{1-\frac{\alpha(v)-u}{\beta(v)-u}
\cdot e^{(\alpha(v)-\beta(v))bqt}}.
\end{eqnarray*}
Hence,
\begin{eqnarray*}
G(t,1,v)=\beta(v)+\frac{\alpha(v)-\beta(v)}{1+\frac{\alpha(v)-1}{1-\beta(v)}\cdot e^{(\alpha(v)-\beta(v))bqt}}.
\end{eqnarray*}
Therefore, by Theorem~\ref{th3.1},
\begin{eqnarray*}
\begin{cases}
g_{0}(t)=G(t,1,0)=\frac{1}{q+p\cdot e^{bt}}\\
g_{n}(t)
=\frac{e^{bt}}{(q+p\cdot e^{bt})^2}\cdot
\int_0^t(q+p\cdot e^{bs})^2e^{-bs}F_n(s)ds, \ n\geq1
\end{cases}
\end{eqnarray*}
with
$$
F_n(t)=bp\delta_{1,n}+bq\sum_{k=1}^{n-1}g_k(t)g_{n-k}(t),\ \ \ n\geq 1.
$$
(i) is proved. (ii) is the direct consequence of (i). The proof is complete. \hfill $\Box$
\end{proof}
\par
\begin{corollary}\label{cor3.4}
Suppose that $X(t)$ is a Markov branching process with $b_0=pb,\ b_3=qb$ $($where $b>0,\ p\in(0,1), p+q=1$$)$, $X(0)=1$. $Y(t)$ is the death number of $\{X(t);t\geq 0\}$ until time $t$. Then the probability generating function of $Y(t)$ is given by
\begin{eqnarray*}
E[v^{Y(t)}]=G(t,1,v)=\sum_{n=0}^{\infty}g_n(t)v^n,
\end{eqnarray*}
where
\begin{eqnarray*}
\begin{cases}
g_{0}(t)=G(t,1,0)=(q+pe^{2bt})^{-1/2}\\
g_{n}(t)
=e^{2bt}\cdot (q+p\cdot e^{2bt})^{-3/2}\cdot
\int_0^te^{-2bs}(q+p\cdot e^{2bs})^{3/2}F_n(s)ds, \ n\geq1
\end{cases}
\end{eqnarray*}
with
$$
F_n(t)=bp\delta_{1,n}+bq\cdot\sum_{k_1,k_2,k_3<n,\ k_1+k_2+k_3=n}g_{k_1}(t)g_{k_2}(t)g_{k_3}(t).
$$
Hence, the probability generating function of death number conditioned on $\tau<\infty$ is given by
\begin{eqnarray*}
E[v^{Y(\tau)}|\tau<\infty]=\beta(v)=\sum_{n=0}^{\infty}g_nv^n,
\end{eqnarray*}
where
\begin{eqnarray*}
\begin{cases}
g_{0}=0,\ \ g_1=p,\\
g_n
=q\cdot \sum_{i+j+k=n}g_{i}g_{j}g_{k}, \ n\geq 2.
\end{cases}
\end{eqnarray*}
\end{corollary}
\par
\begin{proof} Note that $ D=\{0\}$, $B(y)=b(p-y+qy^3)$ and
\begin{eqnarray*}
B_{ D}(y,v)+\bar{B}_{ D}(y)=b(pv-y+qy^3).
\end{eqnarray*}
Let $y(t)=G(t,1,v)=\sum_{n=0}^{\infty}g_n(t)v^n$ be the solution of (\ref{eq3.6}) with $u=1$, then
\begin{eqnarray*}
  \begin{cases}g'_0(t)=b(-g_0(t)+qg^3_0(t)),\\
  g_0(0)=1.
  \end{cases}
\end{eqnarray*}
Solving the above equation yields
\begin{eqnarray*}
g_0(t)=(q+pe^{2bt})^{-1/2}.
\end{eqnarray*}
Therefore, the first result follows from Theorem~\ref{th3.1} and taking limit yields the second result.
The proof is complete. \hfill $\Box$
\end{proof}

\section*{Acknowledgement}
\par
 The work of Yanyun Li is supported by the National Natural Science Foundation of China (NSFC 11931018) and the Guangdong Province Key Laboratory of Computational Science at the Sun Yat-Sen University (2020B1212060032). The work of Junping Li is supported by the National Natural Science Foundation of China (No. 11771452, No. 11971486).


\end{document}